\documentclass{article}
\usepackage{amsmath, amsxtra,amssymb, latexsym, amsfonts, amscd} 
\usepackage{multicol,color}
\usepackage{bbm}
\usepackage{float}
\usepackage{soul}
\usepackage{graphicx}
\usepackage{subcaption}
\definecolor{labelkey}{rgb}{0,0.08,0.45}
\definecolor{refkey}{rgb}{0,0.6,0.0}
\definecolor{Brown}{rgb}{0.45,0.0,0.05}
\definecolor{lime}{rgb}{0.00,0.8,0.0}
\definecolor{lblue}{rgb}{0.5,0.5,0.99}
\usepackage{pstricks-add}
\usepackage{amssymb}
\usepackage{theorem}
\usepackage{fancyhdr}
\usepackage{tikz}
\usepackage{enumerate}
\usepackage[margin=1.1in]{geometry}
\usetikzlibrary{arrows}
\usepackage{hyperref}
\usepackage{esint}
\usepackage{pdfcomment}
\definecolor{labelkey}{rgb}{0.6,0.6,0.6}
\definecolor{refkey}{rgb}{0,0.6,0.0}
\def\disp{\displaystyle}
\def\e{\epsilon}

\def\vt{\vartheta}
\def\dd{\delta}

\def\lm{\lambda}
\def\O{\Omega}
\def\Tilde{\widetilde}
\def\tilde{\widetilde}

\def\({\left(}
\def\){\right)}
\def\[{\left[}
\def\]{\right]}
\def\n{\left \|}
\def\en{\right \|}
\def\nn{\left \{ }
\def\hnn{\right \}}

\def\ox{\bar{x}}
\def\oy{\bar{y}}
\def\oz{\bar{z}}
\def\ov{\bar{v}}
\def\ou{\bar{u}}
\def\ot{\bar{t}}
\def\oT{\overline{T}}

\def\cone{\hbox{}}
\def\span{\hbox{\rm span}}
\def\gph{\hbox{}}

\def\gg{\gamma}
\def\dn{\downarrow}

\def\tto{\rightrightarrows}
\def\st{\stackrel}

\def\Limsup{\mathop{{\rm Lim}\,{\rm sup}}}

\def\hat{\widehat}

\def\tilde{\widetilde}
\def\Tilde{\widetilde}
\def\Bar{\overline}
\def\la{\langle}
\def\ra{\rangle}

\def\B{I\!\!B}
\def\h{\hfill\Box}
\def\R{\mathbb{R}}
\def\N{\mathbb{N}}

\def\co{\mbox{\rm co}\,}
\def\gph{\mbox{\rm gph}\,}
\def\epi{\mbox{\rm epi}\,}

\def\dom{\mbox{\rm dom}\,}

\def\proj{\mbox{\rm proj}\,}

\def\cone{\mbox{\rm cone}\,}

\def\dn{\downarrow}
\def\O{\Omega}

\def\vph{\varphi}
\def\emp{\emptyset}
\def\st{\stackrel}
\def\oR{\Bar{\R}}
\def\lm{\lambda}

\def\gg{\gamma}

\def\dd{\delta}

\def\al{\alpha}
\def\vth{\vartheta}

\def\be{\beta}
\def\vth{\vartheta}
\def\ph{\varphi}

\def\N{I\!\!N}
\def\th{\theta}

\newtheorem{theorem}{Theorem}[section]

\newtheorem{proposition}[theorem]{Proposition}
\newtheorem{definition}[theorem]{Definition}
\theoremstyle{plain}{\theorembodyfont{\rmfamily}
}
\theoremstyle{plain}{\theorembodyfont{\rmfamily}
}
\theoremstyle{plain}{\theorembodyfont{\rmfamily}
}
\theoremstyle{plain}{\theorembodyfont{\rmfamily}
\newtheorem{example}[theorem]{Example}}

\theoremstyle{plain}{\theorembodyfont{\rmfamily}
\newtheorem{remark}[theorem]{Remark}}

\def\eq{\begin{equation}}
\def\eeq{\end{equation}}

\begin{document}
\begin{center}
\end{center}
\begin{center}
{\bf DISCRETE APPROXIMATIONS AND OPTIMALITY CONDITIONS\\ FOR CONTROLLED FREE-TIME SWEEPING PROCESSES}\\[3ex]
GIOVANNI COLOMBO\footnote{Dipartimento di Matematica ``Tullio Levi-Civita", Universit$\grave{\textrm{a}}$ di Padova, 35121 Padova, Italy (colombo@math.unipd.it), and G.N.A.M.P.A. of Istituto Nazionale di Alta Matematica ``Francesco Severi", Piazzale Aldo Moro 5, 00185 -- Roma, Italy. The research of the first author was partly supported by the project funded by the EuropeanUnion – NextGenerationEU under the National
Recovery and Resilience Plan (NRRP), Mission 4 Component 2 Investment 1.1 -
Call PRIN 2022 No. 104 of February 2, 2022 of Italian Ministry of
University and Research; Project 2022238YY5 (subject area: PE - Physical
Sciences and Engineering) ``Optimal control problems: analysis,
approximation and applications''.}
\quad BORIS S. MORDUKHOVICH\footnote{Department of Mathematics, Wayne State University, Detroit, Michigan 48202, USA (boris@math.wayne.edu). Research of this author was partially supported by the USA National Science Foundation under grants DMS-1808978 and DMS-2204519, by the Australian Research Council under grant DP-190100555, and by Project~111 of China under grant D21024.} 
\quad DAO NGUYEN\footnote{Department of Mathematics and Statistics, San Diego State University, San Diego, CA 92182, USA (dnguyen28@sdsu.edu). Research of this author was supported by the AMS-Simon Foundation.}
\quad TRANG NGUYEN\footnote{Department of Mathematics, Wayne State University, Detroit, Michigan 48202, USA (daitrang.nguyen@wayne.edu). Research of this author was partially supported by the USA National Science Foundation under grant DMS-1808978 and DMS-2204519.}
\end{center}
\small{\bf Abstract.} The paper is devoted to the study of a new class of optimal control problems governed by discontinuous constrained differential inclusions of the sweeping type with involving the duration of the dynamic process into optimization. We develop a novel version of the method of discrete approximations of its own qualitative and numerical values with establishing its well-posedness and strong convergence to optimal solutions of the controlled sweeping process. Using advanced tools of first-order and second-order variational analysis and generalized differentiation allows us to derive new necessary conditions for optimal solutions of the discrete-time problems and then, by passing to the limit in the discretization procedure, for designated local minimizers in the original problem of sweeping optimal control. The obtained results are illustrated by a numerical example.\\[1ex]
\vspace{-0.2in}

\section{Introduction, Problem Formulation, and Discussions}\label{intro}
\setcounter{equation}{0}\vspace*{-0.1in}

The {\em sweeping process} (``processus de rafle'') was introduced by Moreau in the 1970s (see \cite{moreau} with more references) and extensively studied by himself and other researchers in the form of   differential inclusion
\begin{equation*}
\left\{\begin{matrix}
\dot{x}(t)\in-N\big(x(t);C(t)\big)\;\textrm{ a.e. }\;t\in[0,T],\\x(0)=x_0\in C(t)\subset\R^n.
\end{matrix}\right.
\end{equation*}
The sweeping process, along with its various modifications, has been utilized in a wide range of fields such that aerospace, process control, robotics, bioengineering, chemistry, biology, economics, finance, management science, and engineering. Additionally, the sweeping process plays a significant role in the theory of variational inequalities and complementarity problems. In the context of mechanical and electrical engineering, this process has been applied to various areas including mechanical impact, Coulomb friction, diodes and transistors as well as queues and resource limits among others. Several notable references for these applications include \cite{aht,mb,brog,BrS,hb1,KP,MM,St}.\vspace*{-0.05in}

Optimal control problems for various types of sweeping processes have been formulated much more recently (see \cite{chhm} and the references therein), while being recognized as very challenging in control theory due to the high discontinuity of the controlled sweeping dynamics and the unavoidable presence of hard state and irregular mixed constraints. Nevertheless, for a rather short period of time, many important results have been obtained on necessary optimality conditions for controlled sweeping processes with valuable applications to friction and plasticity, robotics, traffic equilibria, hysteresis, economics, and other fields of engineering and applied sciences; see, e.g., \cite{adam,ac,bk,b1,CaoMordukhovich2018DCDS,cg21,chhm,cmn18b,cmn18a,pfs,pfs23,her-pall,khalil,mn20,vera,zeidan} with more references and discussions.\vspace*{-0.05in}

Nevertheless, there are great many unsolved problems in optimal control theory for sweeping processes with strong requirements for further applications. Some of these issues are considered in our paper. Specifically, we address here, for the {\em first time in the literature} on controlled sweeping processes, the case where the {\em duration of the dynamical process} is included into optimization and then develop for such problems effective techniques of deriving necessary optimality conditions based on a {\em novel version} of the {\em method of discrete approximations}. This approach, being married to advanced tools of first-order and second-order variational analysis and generalized differentiation, allows us to obtain new necessary conditions for local optimal solutions to both discrete-time and continuous-time frameworks of nonsmooth sweeping dynamics in the presence of general constraints on the trajectory endpoints and process duration. In this way, we solve a {\em long-standing question} in the theory of discrete approximations of optimal control problems (not only of the sweeping type) about handling general {\em measurable} (not just piecewise continuous or of bounded variation) control functions.
Establishing the well-posedness and strong convergence for the novel type of discrete approximations and then passing to the limit from the derived descriptions of their optimal solutions, we arrive at the set of necessary optimality conditions for local minimizers of the original sweeping control problem containing {\em significant new features} in comparison with previously known results; see below and more discussions in Remark~\ref{rem-compar}. The obtained optimality conditions are illustrated by a numerical example of its own interest. Applications of the established results to solving some practical problems from marine surface vehicle and nanoparticle modeling are given in a separate forthcoming paper.\vspace*{-0.03in}

Let us now present a precise mathematical formulation of the problem under consideration in this paper. By $(P)$, we denote the optimal control problem of minimizing the Mayer-type cost functional
\begin{equation}\label{cost1}
J[x,u,T]:=\varphi(x(T),T)
\end{equation}
over control actions $u(\cdot)$ and the corresponding sweeping trajectories $x(\cdot)$ defined on the variable time interval $[0,T]$ and satisfying the constraints
\begin{equation}\label{Problem}
\left\{\begin{matrix}
\dot{x}(t)\in-N\big(x(t);C(t)\big)+g\big(x(t),u(t)\big)\;\textrm{ a.e. }\;t\in[0,T],\;x(0)=x_0\in C(0)\subset\R^n,\\
u(t)\in U\subset\R^d\;\textrm{ a.e. }\;t\in[0,T],\\ (x(T), T)\in \Omega_x\times \Omega_T\subset\R^n\times [0, \infty ), 
\end{matrix}\right.
\end{equation}
where $\O_x$ and $\O_T$ are subsets of $\R^n$ and $[0,\infty)$, respectively, and where $C(t)$ is a convex {\em polyhedron} given by
\begin{equation}\label{C}
\left\{\begin{matrix}
C(t):=\bigcap_{j=1}^{s}C^j(t)\textrm{ with }C^j(t):=\nn x\in\R^n\;\big|\;\la x^j_*(t),x\ra\le c_j(t)\hnn,
\\ 
\|x_\ast^j(t)\|=1$, $j=1,\ldots,s,\,t\in [0,\infty).
\end{matrix}\right.
\end{equation}
Recall that the {\em normal cone of convex analysis} $N(x;C)$ is defined by
\begin{equation}\label{nc}
N(x;C):=\big\{v\in\R^n\;\big|\;\la v,y-x\ra\le 0,\;y\in C\big\}\textrm{ if }x\in C\textrm{ and }N(x;C):=\emp\textrm{ if }x\notin C.
\end{equation}
The latter tells us that problem \eqref{Problem} automatically contains the {\em pointwise state constraints}
\begin{equation}\label{e:8}
x(t)\in C(t),\textrm{ i.e., }\la x^j_*(t),x(t)\ra\le c_j(t)\;\textrm{ for all }\;t\in [0,T]\;\textrm{ (with different } T)\;\textrm{ and }\;j=1,\ldots,s.
\end{equation}
In fact, the sweeping dynamics intrinsically induces {\em irregular mixed constraints} on controls and trajectories that are the most challenging and largely underinvestigated in control theory even in particular settings.\vspace*{-0.05in}

In what follows, we identify the trajectory $x:[0,T]\to\R^n$ with its extension to the interval $(0,\infty)$ given by
\begin{equation*}
x_e(t):=x(T)\;\textrm{ for all }\;t>T
\end{equation*}
and for $x\in W^{1,2}([0,T],\R^n)$ define the norm 
\begin{equation*}
\n x\en_{W^{1,2}}:=\n x(0)\en+\n \dot{x}_e\en_{L_2}.
\end{equation*}
Let us specify what we mean by a $W^{1,2}\times L^2\times \R_+$-local minimizer of (P) and its relaxation; cf. \cite{cmn18a} when the duration of the process is fixed.\vspace*{-0.05in}

\begin{definition}\label{Def3.1}
A feasible solution $(\ox(\cdot),\ou(\cdot),\oT)$ to $(P)$ is a {\sc $W^{1,2}\times L^2\times \R_+$-local minimizer} for this problem if there is $\e>0$ such that $J[\ox,\ou,\oT]\le J[x,u,T]$ for any feasible $(x(\cdot),u(\cdot),T)$ satisfying the constraints in $(P)$ and
\begin{equation*}
\int_0^{\oT}\(\n\dot{x}_e(t)-\dot{\ox}_e(t)\en^2+\n u(t)-\ou(t)\en^2\)dt+(\oT-T)^2<\e.
\end{equation*}
\end{definition}
The {\em relaxed} version $(R)$ of problem $(P)$ is defined as follows: 
\begin{equation*}
\textrm{ minimize }\;\hat{J}[x,u,T]:=\varphi(x(T),T)
\end{equation*}
on absolutely continuous trajectories of the {\em convexified} differential inclusion
\begin{equation}\label{conv}
\dot{x}(t)\in-N\big(x(t);C(t)\big)+\co g\big(x(t),U\big)\;\textrm{ a.e. }\;t\in[0,T],\;x(0)=x_0\in C(t)\subset\R^n,
\end{equation}
where ``co" signifies the convex hull of the set in question.\vspace*{-0.05in}

\begin{definition}\label{relaxed} Let $(\ox(\cdot),\ou(\cdot),\oT)$ be a feasible solutions to $(P)$. We say that it is a
{\sc relaxed $W^{1,2}\times L^2\times \R_+$-local minimizer} for $(P)$ if there is $\e>0$ such that
\begin{equation*}
\ph\big(\ox(\cdot),\oT\big)\le\ph\big(x(\cdot),T\big)\;\textrm{ whenever }\;\int_0^{\oT}\(\n\dot{x}_e(t)-\dot{\ox}_e(t)\en^2+\n u(t)-\ou(t)\en^2\)dt+(\oT-T)^2<\e,
\end{equation*}
where $u(t)\in\co U$ a.e.\ on $[0,T]$ with $u(\cdot)$ being a measurable control function, and where $x(\cdot)$ is a corresponding relaxed trajectory of the convexified inclusion \eqref{conv} that can be uniformly approximated in $W^{1,2}([0,T];\R^n)$ by feasible trajectories to $(P)$ generated by piecewise constant controls $u^k (\cdot)$ on $[0,T]$ whose convex combinations converges to $u(\cdot)$ in the norm topology $L^2([0,T];\R^d)$.
\end{definition}\vspace*{-0.1in}
Of course, there is no difference between $W^{1,2}\times L^2\times\R_+$-local minimizers for $(P)$ and its relaxed counterpart $(R)$ in the case of convex problems, but it is also the case for a broad range of problems without any convexity. We refer the reader to \cite{cmn18a} for more details and discussions on this topic.\vspace*{-0.05in}

Observe that when $C(t)\equiv\R^n$, problem $(P)$ reduces to the standard framework of {\em nonsmooth optimal control} of ODE systems since our standing assumptions formulated below impose the Lipschitz continuity of $g$ in $x$. When $g\equiv 0$, problem $(P)$ takes the form of {\em optimal control of differential inclusions} $\dot x(t)\in F(t,x(t))$ with $F(x,t):=-N(x;C(t))$. However, the principal difference between the sweeping control case of $(P)$ and the well-developed control theory for differential inclusions (see, e.g., \cite{cl,Mord,v} and the reference therein) is that the latter theory crucially depends on the {\em Lipschitz continuity} of the multifunction $x\mapsto F(x,t)$, which is {\em never the case} of the normal cone mapping in $(P)$ generated by a nontrivial (constant or moving) set $C(t)$.\vspace*{-0.05in}

Our approach to investigate problem $(P)$ with establishing necessary optimality conditions for its local minimizers is based on the {\em method of discrete approximations} developed in \cite{m95,m96,Mord} for Lipschitzian differential inclusions and then extended in \cite{b1,CaoMordukhovich2018DCDS,chhm,cmn18b,cmn18a,pfs} to various kinds of controlled sweeping processes. Note that \cite{m96} is the only publication employing discrete approximations to handle free-time problems while in the case of Lipschitzian inclusions. Our method consists of constructing well-posed discrete approximations of $(P)$ whose solutions strongly converge to the prescribed local minimizers of $(P)$, then deriving necessary optimality conditions for discrete-time problems, and finally establishing, by passing to the limit with the diminishing step of discretization, necessary optimality conditions for local minimizers of $(P)$.\vspace*{-0.03in}

The paper is organized as follows. Below in this section, we formulate assumptions on the problem data and present a recent theorem on the existence of optimal solution to $(P)$ that supports our analysis.\vspace*{-0.05in}

Section~\ref{Dis-app} is devoted to the construction and justification of {\em well-posed discrete approximations} of the continuous-time problem $(P)$. The first novel result here of its own interest establishes a {\em strong $W^{1,2}\times L^2$-approximation} of any feasible solution $(x,u)\in W^{1,2}([0,T];\R^n)\times L^2([0,T];\R^d)$ to $(P)$ by appropriately extended solutions to the discretized sweeping process. Based on this crucial result, we construct a well-posed sequences of discrete optimal control problems whose optimal solutions strongly converge to the designated local minimizer of $(P)$.\vspace*{-0.05in}

It occurs that even if the initial data of $(P)$ are differentiable or/and convex, the discrete-time approximating problems inevitably become {\em nonsmooth} and {\em nonconvex} due to geometric constraints generated by the {\em graph} of the normal cone \eqref{nc}. To tackle these challenges, we need appropriate tools of {\em variational analysis} including {\em first-order} and {\em second-order generalized differentiation}. Section~\ref{Var-Ana} overviews the required constructions and results.\vspace*{-0.05in}

In Section~\ref{sec:NCO}, we employ variational analysis and generalized differentiation to derive {\em necessary optimality conditions} for {\em discrete-time} problems. This is done by reducing the problems with discrete dynamics to nondynamical problems of {\em nondifferentiable programming} with increasingly many geometric constraints by employing {\em full calculus} of generalized differentiation taken from \cite{m-book2} and explicit second-order calculations associated with the sweeping dynamics. Due to the strong convergence of discrete optimal solutions established in Section~\ref{Dis-app}, the results obtained here can be viewed as {\em suboptimality conditions} for the original sweeping control problem.\vspace*{-0.05in} 

The culmination of our study is Section~\ref{nc-sweep}, where we utilize the stability of discrete approximations and the robustness of generalized differential constructions for {\em passing to the limit} from the necessary optimality conditions for discrete-time problems. Establishing an appropriate convergence of {\em adjoint arcs} occurs to be the most challenging aspect of this process. As a result, we obtain {\em new necessary conditions expressed in terms of the original problem data} for the designated local minimizers of $(P)$. Due to the characteristic features of the sweeping dynamics, the resulting optimality conditions include {\em signed measures}, not just nonnegative ones, which may present significant difficulties for implementation. Nevertheless, the newly obtained {\em support condition} allows us to largely overcomes this obstacle. This condition is obtained by exploiting the particular structure of the approximate adjoint equation, together with a precise calculation of the coderivative of the normal cone mapping.
Section~\ref{examples} contains a numerical example showing how the obtained optimality conditions make it possible to completely solve to find the original sweeping control problem. The paper is finished with Section~\ref{sec:Conclusions} containing concluding remarks and discussions of some topics of our future research.\vspace*{-0.03in}

In the reminder of this section, we formulate and discuss {\em standing assumptions} for the rest of the paper. Consider the collection of {\em active indices} of inequality constraints defined by
\begin{equation*}
I(t,x):=\{j\in\{1,\ldots,s\}\;|\;\la x^j_{\ast}(t),x\ra=c_j(t)\}.
\end{equation*}\vspace*{0.03in}
{\bf(H1}) The set $U\ne\emp$ is closed and bounded in $\R^d$. The generating functions $x^j_{\ast}(\cdot)$, and $c_j(\cdot)$ are Lipschitz continuous with a common Lipschitz constant $L$.\\\vspace*{0.03in}
{\bf(H2)} The {\em uniform Slater condition} is satisfied:
\begin{equation}\label{sla}
\mbox{for every }\;t\in[0,T]\;\mbox{  there exists }\,x\in\R^n\;\mbox{ such that }\;\la x^j_{\ast}(t), x\ra<c_j(t)\;\mbox{ whenever }\;i=1,\ldots, s.
\end{equation}
This condition yields the {\em positive linear independence constraint qualification} (PLICQ) along $x(t)$ on $[0,T]$ with a varying time $T$ formulated as follows:
\begin{equation}\label{plicq}
\Big[\sum_{j\in I(t,x)}\lambda_jx^j_*(t)=0,\;\lambda_j\in\R_+\Big]\Longrightarrow\big[\lambda_j=0\;\textrm{ for all }\;j\in I(t,x)\big].
\end{equation}
In \cite{CCMN21,hjm}, the reader can find more discussions on this and related topics. Recall that the (stronger) {\em linear independence constraint qualification} (LICQ) holds if the restriction $\lm_j\in\R_+$ in \eqref{plicq} is dropped.\\\vspace*{0.03in}
{\bf(H3)} The perturbation mapping $g\colon\R^n\times U\to\R^n$ is {\em Lipschitz continuous} with respect to $(x,u)$ whenever $u\in U$ and $x$ belongs to a bounded subset of $\R^n$ satisfying there the {\em sublinear growth condition}
\begin{equation*}
\|g(x,u)\|\le\be\big(1+\|x\|\big)\;\mbox{ for all }\;u\in U
\end{equation*}
with some given positive constant $\be$.\vspace*{-0.03in}

Define further set-valued mapping $F\colon[0,T]\times\R^n\times\R^d\tto\R^n$ by
\begin{equation}\label{F0}
F(t,x,u):=N(x;C(t))-g(x,u)
\end{equation}
and deduce from the classical Motzkin theorem of the alternative the representation
\begin{equation}\label{F}
F(t,x,u)= \Big\{\sum_{j\in I(t,x)}\lm^j x^j_*(t)\;\Big|\;\lm^j\ge 0\Big\}-g(x,u),
\end{equation}
where the moving set $C(t)$ is taken from \eqref{C}.
Then the sweeping differential inclusion \eqref{Problem} can be rewritten as 
\begin{equation}
\label{e:SP1}
\(-\dot x(t),u(t)\) \in F(t,x(t),u(t)) \times U \;\;\mbox{a.e.}\;\; t\in [0,T].
\end{equation}\vspace{-0.15in}

The following theorem is taken from \cite[Theorem~3.1]{hjm}.\vspace*{-0.1in}
\begin{theorem}\label{ext}
Let $(P)$ be the sweeping optimal control problem with the equivalent form \eqref{e:SP1} of the sweeping differential inclusion over all the $W^{1,2}([0,T];\R^n) \times L^2([0,T];\R^d)\times [0,\infty)$ triples $(x(\cdot),u(\cdot),T)$, and let the moving polyhedron satisfy the uniform Slater condition \eqref{sla}. Then (P) admits a Lipschitzian optimal solution $x(\cdot)$.  
\end{theorem}\vspace*{-0.25in}

\section{Well-Posed Discrete Approximations}\label{Dis-app}\vspace*{-0.05in}

In this section, we aim at developing well-posed {\em discrete approximations} of the sweeping control problem $(P)$, which deals with constrained differential inclusions involving free time. This method, providing finite-dimensional approximations of infinite-dimensional continuous-time problems, definitely has some numerical flavor, while our main attention here is paid to using the discrete approximation approach to derive efficient necessary optimality conditions for the original problem $(P)$.\vspace*{-0.05in} 

To simplify the exposition, we employ the {\em explicit Euler scheme} to replace the time derivative $\dot x(t)$ in \eqref{Problem} by the sequence of finite differences
\begin{equation*}
\dot x(t)\approx\frac{x(t+h)-x(t)}{h}\;\textrm{ as }\;h\dn 0,
\end{equation*}
which is formalized as follows: for each $k\in\N$, consider a  real number $T_k$ approximating $T$ and the uniform grid
\begin{equation}\label{grid}
\left\{\begin{matrix}
t^k_0:=0,\quad t^k_k:=T_k,\\
t^k_{i+1}:=t^k_i+h^k,\;i=0,1,\ldots,k-1,
\end{matrix}\right.
\end{equation}
with $h^k:=T_k/k$ and $\N:=\{1,2,\ldots\}$ standing for the collection of natural numbers.\vspace*{-0.05in}

The following theorem plays a principal role in the subsequent developments. It justifies the possibility to {\em $W^{1,2}\times L^2$-strongly} approximate {\em any feasible solution} to $(P)$ by a sequence of extended solutions to certain {\em discretized perturbed} sweeping processes. To proceed, recall some notation. Given a Lebesgue measurable set $S\subset\R^n$ with  positive finite measure $|S|$ and a Lebesgue measurable mapping $f\colon\R^n\to\R^m$, denote its {\em average} by
\begin{equation*}
\fint_S f(s):=\frac{1}{|S|}\int_S f(s)ds.
\end{equation*}
Further, having some sets $\Omega,S\subset\R^n$ with $\Omega\subset S$, define the characteristic function $\mathbbm{1}_{\Omega}(x):=1$ if $x\in\Omega$ and $\mathbbm{1}_{\O}(x):=0$ if $x\in S\setminus\O$. \vspace*{-0.05in}

\begin{theorem}\label{Thr1} Let $(\ox(\cdot),\ou(\cdot))\in W^{1,2}([0,\oT];\R^n)\times L^2([0,\oT];\R^d)$ be a feasible solution to problem $(P)$ under the assumptions in $\rm{(H1)}$ and ${\rm(H3)}$. Then the following assertions hold:\\
{\bf(i)} There exists a sequence of piecewise constant control functions $\{u^k\;|\;k\in\N\}$ defined on $[0,\oT]$ such that $u^k(\cdot)$ converges to $\ou(\cdot)$ strongly
in the $L^2$-norm topology on $[0,\oT]$ as $k\to\infty$. \\
{\bf(ii)}) There exists a sequence of piecewise linear functions $\{x^k\;|\; k\in \N\}$ converging strongly to $\ox(\cdot)$
in the $W^{1,2}$-norm topology on $[0,\oT]$ and such that $x^k(0)=\ox(0)$ for all $k\in \N$ and all $i=0,\ldots,k-1$ while satisfying
\begin{equation}
\label{e:x-dc}
\dot x^k(t)\in -\mbox{N}(x^k(t^k_i);C^k_i)+g(x^k(t^k_i), u^k(t^k_i))+\tau^k_i\B, \;t\in[t^k_i,\;t^k_{i+1}),
\end{equation}
where $\tau^k_i\geq 0$ with $h^k\sum^{k-1}_{i=1}\tau^k_i\to 0$ as $k\to\infty$, and where the perturbed polyhedra $C^k_i$ are defined by
\begin{equation}\label{C_i}
C^k_i:=\disp\bigcap_{j=1}^s C^k_{ij}\;\mbox{ with }\;C^k_{ij}:=\big\{x\in\R^n\;\big|\;\la x,x^j_{\ast}(t^k_i)\ra \leq c^{ij}_k\big\}
\end{equation}
being generated by the vectors $x^j_{\ast}(t^k_i)$ and suitable numbers $c^{ij}_k$.\\
{\bf(iii)} The piecewise linear extensions of $c^{ij}_k$ and $x^j_{\ast}(t^k_i)$ converge uniformly on $[0,\oT]$ to $c_j(t)$ and $x^j_{\ast}(t)$ from \eqref{C}, respectively. Finally, all $x^k(\cdot)$ are Lipschitz continuous on $[0,\oT]$ with the same Lipschitz constant as $\bar{x}(\cdot)$.
\end{theorem}\vspace*{-0.1in}
\begin{proof}$\,$
Let us first construct control functions $u^k(\cdot)$ to approximate $\ou(\cdot)$. Applying \cite[Lemma~6]{cg21} to $\ou\in L^2([0,\oT];\R^n)$ and considering the control sequence
$$u^k(t):=\disp\sum^{k-1}_{i=0}\fint_{t^k_i}^{t^k_{i+1}}\ou(s)ds \mathbbm{1}_{[t^k_i,t^k_{i+1})}(t),\quad k\in\N,$$
we obtain that $u^k\to\ou$ a.e.\ on $[0,\oT]$ as $k\to\infty$. Define further the sequences
\begin{equation}\label{vel}
v^k_i:=\dfrac{\ox(t^k_{i+1})-\ox(t^k_{i})}{h^{k}}, 
\end{equation}
\begin{equation}\label{tra}
x^k(t):=\ox(t^k_i)+(t-t^k_i)v^k_i= \ox(t^k_i)+(t-t^k_i)\disp \fint_{[t^k_i,t^k_{i+1})}\dot \ox (s)ds,\;t\in [t^k_i,t^k_{i+1}).
\end{equation}
Fixing $i = 0,\ldots,k- 1$, denote $I^i_k:=[t^k_i,t^k_{i+1})$, $\sigma_k(t):=\max\{t^k_i\;|\;i=0,...,k-1,\;t^k_i \leq t\}$, and 
$$
c^{ij}_k:= \left\{\begin{array}{ll}
\disp \max_{t\in [t^k_i,t^k_{i+1}]}c_j(t) &\mbox{if }\;\la \ox(t), x^j_{\ast}(t)\ra < c_j(t)\;\mbox{ for }\;t\in [t^k_i, t^k_{i+1}],\\
\la \ox(t^k_i), x^j_{\ast}(t^k_i)\ra &\mbox{ otherwise} \end{array}\right.
$$
whenever $j = 1,\ldots, s$. For such $j$, consider further the sequence
$$
c^{j}_k(t):=\disp\sum^{k-1}_{i=1}c^{ij}_k \mathbbm{1}_{[t^k_i,t^k_{i+1})}(t).
$$
We intend to show that $c^j_k(t) \to c_j(t)$ uniformly in $[0,\oT]$ as $k\to \infty$. To proceed, for every $t\in [0,\oT]$ define $I_k(t):=[\sigma_k(t),\sigma_k(t)+h^k]$ and then fix
$j\in\{1,\ldots ,s\}$ and $\ot \in [0,\oT]$. If $\la \ox(\ot),x^j_{\ast}(\ot)\ra<c_j(\ot)$, we get  $\la \ox(t),x^j_{\ast}(t)\ra  < c_j(t)  $ for all $t \in [\sigma_k(\ot),\sigma_k(\ot)+ h^k]$ provided that $k$ is sufficiently large. This tells us that $c^j_k(\ot)=\disp \max_{t\in [\sigma_k(\ot),\sigma_k(\ot)+ h^k]}c_j(t)=c_j(\tilde t)$ for some $\tilde t=\tilde t(k,j,\ot)\in [\sigma_k(\ot),\sigma_k(\ot)+ h^k]$, and hence 
$$|c_j(\ot)-c^j_k(\ot)|\leq |c_j(\ot)-c_j(\tilde t)|\leq Lh^k.$$
Letting $\la \ox(\ot),x^j_{\ast}(\ot)\ra  = c_j(\ot)$, it follows that $c^j_k(\ot)=c^j_k(\sigma_k(\ot))$, 
$$
\la \ox(\sigma_k(\ot)),x^j_{\ast}(\sigma_k(\ot))\ra = c^j_k(\sigma_k(\ot)),\;\mbox{ and }\;\la \ox(\ot), x^j_{\ast}(\ot
)\ra = c_j(\ot),$$
which brings us to the relationships
\begin{equation*}
\begin{array}{ll}
&|c_j(\ot)-c^j_k(\sigma_k(\ot))|= |\la \ox(\ot),x^j_{\ast}(\ot)\ra -\la \ox(\sigma_k(\ot)), x^j_{\ast}(\sigma_k(\ot))\ra |\le|\la \ox(\ot),x^j_{\ast}(\ot)- x^j_{\ast}(\sigma_k(\ot))\ra\\
&+\la \ox(\ot)-\ox(\sigma_k(\ot)), x^j_{\ast}(\sigma_k(\ot))\ra|\leq Kh^k
\end{array}
\end{equation*}
for a suitable constant $K$ independent of $k$ and $\ot$. 
Define $g_k(t):= g(\ox(\sigma_k(t)), u^k(t))$ and $\zeta_k(t):= g_k(t)-w_k(t)$, where $ w_k(t):=\displaystyle \sum^{k-1}_{i=0}v^k_i\mathbbm{1}_{I^i_k}(t)$, $t\in [0,\oT]$. 
It follows from the normal cone representation in \eqref{Problem} that, a.e. on $[0,\oT]$,
$$g(\ox(t), \ou(t))-\dot\ox(t)=\sum _{j=1}^s\lm_j(t)x^j_{\ast}(t),$$ 
where $\lm_j(t)\geq 0$ and $\lm_j(t)=0$ for each $t$ with $\la \ox(t), x^j_{\ast}(t)\ra <c_j(t)$. Moreover, it is easy to check that the numbers $\|\lm_j\|_{L^{\infty}}$ are finite since $\ox$ is Lipschitz continuous, $g$ is bounded, and all the active inequality constraints are positively linearly independent. Observe that
\begin{equation}\label{2.16bis}
\begin{aligned}
\zeta_k(t) &= \fint_{[\sigma_k(t),\sigma_k(t)+ h^k)}[g(\ox(r),\ou(r))+\dot \ox(r)]dr - \fint_{[\sigma_k(t),\sigma_k(t)+ h^k)}[g(\ox(r),\ou(r))-g_k(r)]dr \\
&= \sum _{j=1}^s \fint_{[\sigma_k(t),\sigma_k(t)+ h^k)} x^j_{\ast}(r)\lm_j(r)dr -\fint_{[\sigma_k(t),\sigma_k(t)+ h^k)}[g(\ox(r),\ou(r))-g_k(r)]dr
\end{aligned}
\end{equation}
and define further the sequence of moving sets
\begin{equation*}
C^k_j(t):= \disp \sum^{k-1}_{i=0} C^k_{ij}\mathbbm{1}_{I^i_k}(t),\quad t\in[0,\oT],
\end{equation*}
where $C^k_i$ is taken from \eqref{C_i}.
In this notation, we have the relationships
\begin{equation}\label{2.19}
\begin{aligned}
&\sum _{j=1}^s \fint_{[\sigma_k(t),\sigma_k(t)+ h^k)} x^j_{\ast}(r)\lm_j(r)dr 
\\
&=\disp\sum_{j=1}^s x^j_{\ast}(\sigma_k(t))\fint_{[\sigma_k(t),\sigma_k(t)+ h^k)}\lm_j(r)dr+\sum _{j=1}^s \fint_{[\sigma_k(t),\sigma_k(t)+ h^k)} \lm_j(r)[x^j_{\ast}(r)-x^j_{\ast}(\sigma_k(t))]dr\\
&\in N(x^k(\sigma_k(t));C^k_j(\sigma_k(t)))+\tilde Lh^k\B= N(x^k(t);C^k_j(t))+\tilde Lh^k\B
\end{aligned}
\end{equation}
for a.e. $t\in[0,\oT]$ and a suitable constant $\tilde L$. On the other hand, using the fact that $\sigma_k(s)=\sigma_k(t)$ for all $s\in I_k(t)$ and $t\in [0,\oT]$ leads us to the equality
$$
\begin{aligned}
\fint_{[\sigma_k(t),\sigma_k(t)+ h^k)}[g(\ox(r),\ou(r))-g_k(r)]dr
&= \fint_{[\sigma_k(t),\sigma_k(t)+ h^k)}[g(\ox(r),\ou(r))-g(\ox(\sigma_k(t), \ou(r)))]dr\\
&+\fint_{[\sigma_k(t),\sigma_k(t)+ h^k)}[g(\ox(\sigma_k(t)), \ou(r))-g_k(r)]dr.
\end{aligned}
$$
The uniform continuity of g with respect to $x$ and $u$ and the convergence of $u^k$ to $\ou$ tell us that the latter expression tends to zero as $k \rightarrow\infty$. Combining this with \eqref{2.16bis} and \eqref{2.19}, we find a sequence $\{\tau^k_i\}$ such that $\tau^k_i\downarrow  0$ as $k \rightarrow\infty$ and \eqref{e:x-dc} is satisfied. The remaining part of the statement is an immediate consequence of the constructions.
\end{proof}

Having in hand the approximation results of Theorem~\ref{Thr1} allows us to build a sequence of {\em discrete-time} optimal control problems, which admit optimal solutions that {\em strongly converge} to a {\em designated local minimizer} of the original sweeping control problem $(P)$.\vspace*{-0.05in} 

Let $(\ox(\cdot),\ou(\cdot),\oT)$ be a relaxed $W^{1,2}\times L^2\times\R_+$-local minimizer of problem $(P)$. Given $\e>0$ from Definition~\ref{relaxed} and any fixed $k\in\N$, formulate the {\em discrete-time} optimal control problem $(P_k)$ as follows:
\begin{eqnarray}\label{d_a_p}
\begin{array}{ll}
\mbox{minimize}\;\;J_k[x^k,u^k,T_k]:=\varphi(x^k_k,T_k)+(T_k-\oT)^2\\
\disp+\sum_{i=0}^{k-1}\int_{t^k_i}^{t^k_{i+1}}\left(\n \frac{x^k_{i+1}-x^k_i}{h^k}-\dot{\bar{x}}(t)\en^2+\n u^k_i-\bar{u}(t)\en^2\right)dt
\end{array}
\end{eqnarray}
over $(x^k,u^k,T_k):=(x^k_0,x^k_1,\ldots,x^k_{k-1}, x^k_k,u^k_0,u^k_1,\ldots,u^k_{k-1},T_k)$ satisfying the constraints
\begin{equation}\label{re_1}
x^k_{i+1}- x^k_i\in -h^kF^k_i(t^k_i,x^k_i,u^k_i)\;
\textrm{ for }\;i=0,\ldots,k-1,
\end{equation}
\begin{equation}\label{re_2}
x^k_0:=x_0\in C(0),
\end{equation}
\begin{equation}\label{Omega}
(x^k_k,T_k)\in\O^k_x\times \O^k_T:= (\O_x+\delta^k\B)\times (\O_T+\delta^k),
\end{equation}
\begin{equation}\label{re_3}
\sum_{i=0}^{k-1}\int_{t^k_i}^{t^k_{i+1}}\(\n\frac{x^k_{i+1}-x^k_i}{h^k}-\dot{\ox}(t)\en^2+\n u^k_i-\ou(t)\en^2\)dt\le \e,
\end{equation}
\begin{equation}\label{4.4bis}
\left\| \frac{x^k_{i+1}-x^k_i}{h^k}\right\| \le L +1\;\text{ for }\;i=0,\ldots , k-1,
\end{equation}
\begin{equation}\label{re_4}
u^k_i\in U\;\textrm{ for }\;i=0,\ldots,k-1,
\end{equation}
\begin{equation}\label{re_5}
\n \(x^k_i,u^k_i\)-\(\ox(t^k_i),\ou(t^k_i)\)\en \le \e\;\textrm{ for }\; i=0,\ldots,k-1,
\end{equation}
\begin{equation}\label{2.27}
|T_k-\oT|\le \e,
\end{equation}
\begin{equation}\label{re_6}
\la x^j_\ast(T_k), x^k_k\ra \leq c_j(T_k)\;\textrm{ for  }\;j =1,\ldots,s,
\end{equation}
where $L$ is the Lipschitz constant of $\bar{x}$, $\delta^k:=\|\ox(\oT)-\hat x^k(\oT)\|$,
\begin{equation}\label{F^k_i}
F^k_i(t^k_i,x^k_i,u^k_i):=N(x^k_i;C^k_i)-g(x^k_i, u^k_i)-\tau^k_i\B,
\end{equation}
and the sequence of piecewise linear functions $\{\hat x^k\}$ is generated by the optimal trajectory $\ox(\cdot)$ according to Theorem~{\rm\ref{Thr1}.\vspace*{-0.03in}

To proceed further, let us add to the standing assumptions (H1)--(H3) the following one:\\[1ex]
{\bf (H4)} The set $\O_x\times\O_T$ is closed around $(\ox(\oT),\oT)$.\vspace*{-0.02in}

In our approach to deriving necessary optimality conditions for the original problem $(P)$, we need to be sure that approximating problems $(P_k)$ admits optimal solutions for large $k$. This can be easily deduced from Theorem~\ref{Thr1} and the classical Weierstrass existence theorem.\vspace*{-0.1in}

\begin{proposition}\label{Thr2}
Let the cost function $\varphi$ be lower semicontinuous on $\R^n\times [0,\oT+\e]$ in addition to the assumptions in {\rm(H1)}, {\rm(H3)}, and {\rm(H4)}. Then the discrete approximation problem $(P_k)$ admits an optimal solution $\{(\ox^k(t^k_i)\;|\;i=0,\ldots,k\}$ whenever $k\in\N$ is sufficiently large.
\end{proposition}\vspace*{-0.1in}
\begin{proof} \; It follows from Theorem~\ref{Thr1} that the set of feasible solutions to problem $(P_k)$ is nonempty for any large $k$. It is obvious from the construction of $(P_k)$ that the feasible solution sets are bounded for each $k\in\N$. Furthermore, the constraint structures of $(P_k)$ and the robustness (closed-graph) property of the set-valued mapping \ref{nc} ensure the closedness of the feasible solution sets. The imposed lower semicontinuity of $\ph$ ensures that the cost function in each $(P_k)$ is also lower semicontinuous for each $k\in\N$. Thus the claimed existence result follows from the classical (one-sided) Weierstrass theorem. 
\end{proof}\vspace*{-0.05in}

Now we are ready to obtain a major result establishing the {\em $W^{1,2}\times L^2\times\R_+$-strong convergence} of {\em optimal solutions} for the discrete-time problems $(P_k)$ to the designed local minimizer $(\ox(\cdot),\ou(\cdot),\oT)$ of the sweeping control problem $(P)$. This theorem makes a bridge between solving problem $(P)$ via its discrete finite-dimensional counterparts $(P_k)$ and will be an important ingredient in what follow to derive optimality conditions in $(P)$ by passing to the limit from those for the discrete approximations.\vspace*{-0.1in}

\begin{theorem}\label{Thr3} Let $(\ox(\cdot),\ou(\cdot),\oT)$ be a relaxed $W^{1,2}\times L^2\times\R_+$--local minimizer for the original problem $(P)$, and let $\varphi$ be continuous around $(\ox(\oT),\oT)$ under the notation and assumptions of Proposition~{\rm\ref{Thr2}}. Then for any extended sequence of optimal solutions $(\ox^k(\cdot),\ou^k(\cdot),\oT_k)$ to $(P_k)$, we have the convergence $ \ox^k(\cdot)\to \ox(\cdot)$ in the $W^{1,2}$-norm topology on $[0,\oT]$, $ \ou^k(\cdot)\to\ou(\cdot)$ in the $L^2$-norm topology on $[0,\oT]$, and  $\oT_k\to \oT$ as $k\to \infty$. Moreover, the Lipschitz constant of such optimal solutions tends to $L$ as $k\to\infty$.
\end{theorem}\vspace*{-0.1in}
\begin{proof}\;
We begin with verifying the inequality
\begin{equation}\label{J}
\disp \limsup_{k\to\infty}J_k[\ox^k,\ou^k,\oT_k]\le J[\ox,\ou,\oT]
\end{equation}
for any sequence of optimal solutions to $(P_k)$. Assuming the contrary gives us a sequence of natural numbers $k\to\infty$ and a positive number $\gamma$ such that 
\begin{equation}\label{cont_J}
J[\ox,\ou,\oT]=\varphi(\ox(\oT),\oT)<J_k[\ox^k,\ou^k,\oT_k]-\gamma\;\text{ for all }\;k.
\end{equation}
Let $(x^k,u^k)$ be the sequence of approximate solutions constructed in Theorem~\ref{Thr1}.  Since $x^k_0=\ox_0$ and $\varphi$ is continuous around $(\ox(\oT),\oT)$, we have the convergence
$$\varphi(x^k_k,\oT_k)\to \varphi(\ox(\oT),\oT)\;\textrm{ as }\;k\to 0.
$$
We deduce from the extension rules in \eqref{vel}, \eqref{tra} and the strong convergence in Theorem~\ref{Thr1} that 
\begin{equation*}
\sum_{i=0}^{k-1}\int_{t^k_i}^{t^k_{i+1}}\(\n\frac{x^k_{i+1}-x^k_i}{h^k}-\dot{\ox}(t)\en^2+\n u^k_i-\ou(t)\en^2\)dt =\int^{\oT}_{0}\(\n\dot{x}^k(t)-\dot{\ox}(t)\en^2+\n u^k_i-\ou(t)\en^2\)dt\to 0
\end{equation*}
as $k\to\infty$, which implies in turn that
$$
J_k[x^k(\cdot),u^k(\cdot),\oT]\rightarrow J[\ox,\ou,\oT]\;\;\mbox{as}\;\;k
\rightarrow\infty.
$$
Using the latter and the feasibility of $(x^k(\cdot),u^k(\cdot),T_k)$ for $(P_k)$ tells us that
\eqref{cont_J} contradicts the optimality of $(\ox^k(\cdot), \ou^k(\cdot),\oT)$ for problems $(P_k)$ when $k$ is sufficiently large.\vspace*{-0.05in}

To complete the proof of the theorem, it remains to show that 
\begin{equation}\label{e:0.32}
\lim_{k\rightarrow\infty}\[\zeta_{k}:=|\oT_k
-\oT|^2+\int_{0}^{\oT_k}\(\n\dot{\ox}^k(t)-\dot{\ox}(t)\en^{2}+\n\ou^k(t)-\ou(t)\en^2\)dt\]=0.
\end{equation}
If it were not so, we would consider any limiting point $\zeta>0$ of the sequence $\{\zeta_k\}$ in \eqref{e:0.32}. For the simplicity of exposition, assume that $ \zeta_k\to \zeta$ as $k \to \infty$. It follows from the boundedness of $\{\oT_k\}$ that there exists $\tilde{T}$ such that $\tilde{T}\in \R$, $ \tilde{T}\le \oT+\e$, and $\oT_k\to \tilde{T}$ as $k\to \infty$. Consider the extended discrete trajectories $\ox^k(t)$ and the extended discrete controls $\ou^k(t)$ on $[0,\tilde{T}]$ defined by
$$
\ox^k(t):=\ox^k(T_k)\;\textrm{ and }\;\ou^k(t):=\ou^k(T_k)\; \textrm{ for }\;t\in(T_k,\tilde{T}] \;\textrm{ when }\;T_k < \tilde{T}.
$$
Due to \eqref{re_3} and \eqref{re_4}, the sequence of extended optimal solutions $\{(\dot\ox^k(\cdot),\ou^k(\cdot),\oT_k)\}$ to $(P_k)$ is bounded in the reflexive space $L^2([0,\tilde{T}];\R^n)\times L^2([0,\tilde{T}];\R^d)\times\R_+$, and thus it contains a weakly converging subsequence in this space, without relabeling. Denote by $(\Tilde v(\cdot),\Tilde u(\cdot),\Tilde T)$ the limit of the selected subsequence and then let
\begin{equation*}
\tilde{x}(t):=x_0+\int_0^t\tilde v(\tau)d\tau\;\textrm{ for all }\;t\in[0,\tilde{T}].
\end{equation*}
Since $\dot{\tilde{x}}(t)=\tilde v(t)$ for a.e.\ $t\in[0,T]$, we get without relabeling that
\begin{equation*}
\big(\ox^k(\cdot),\ou^k(\cdot),\oT_k\big)\to\big(\tilde x(\cdot),\tilde u(\cdot),\tilde T\big)\;\textrm{ as }\;k\to\infty
\end{equation*}
in the  weak topology of $W^{1,2}([0,\tilde{T}];\R^n)\times L^2([0,\tilde{T}];\R^d)\times\R_+$. Then Mazur's weak closure theorem gives us a sequence of convex combinations of $(\ox^k(\cdot),\ou^k(\cdot),\oT_k)$, which converges to $(\tilde x(\cdot),\tilde u(\cdot),\tilde T)$ strongly in $W^{1,2}([0,\tilde{T}];\R^n)\times L^2([0,\tilde{T}];\R^d)\times\R_+$, and hence $(\dot{\ox}^k(t),\ou^k(t),\oT_k)\to(\dot{\tilde x}(t),\tilde u(t),\tilde T)$ for a.e. $t\in[0,\tilde{T}]$ along a subsequence. The obtained pointwise convergence of convex combinations allows
us to conclude that $\tilde u(t)\in\co U$ for a.e.\ $t\in[0,\tilde{T}]$ and that $\tilde x(\cdot)$ satisfies the convexified differential inclusion \eqref{conv}.\vspace*{-0.05in}

Using the sign "$\sim$" for expressions which are equivalent
as $k\rightarrow\infty$ brings us to
\begin{equation*}
\begin{array}{ll}
\disp\sum_{i=0}^{k-1}\int_{t^k_i}^{t^k_{i+1}}\(\n\frac{\ox^k_{i+1}-\ox^k_i}{h^i_k}-\dot{\bar{x}}(t)\en^{2}+\n u^k_i-\ou(t)\en^2\)dt &=\disp\int_{0}^{\oT_k}\(\n\dot{\bar{x}}_k
(t)-\dot{\bar{x}}(t)\en^{2}+\n u^k_i-\ou(t)\en^2\)dt\\
&\sim\disp\int_{0}^{\tilde{T}}\(\n\dot{\bar{x}}_k(t)-\dot{\bar{x}}(t)\en^{2}+\n u^k_i-\ou(t)\en^2\)dt
\end{array}
\end{equation*}
as $k\rightarrow\infty$.
Invoking the convexity in $v$ of the function $f(v,t):=\n v-\dot{\ox}(t)\en^2$, we get that the integral functional 
$$I[v]:=\int_{0}^{\tilde{T}}\n v(t)-\dot{\bar{x}}(t)\en^{2}dt$$
is lower semicontinuous in the weak topology of $L^{2}[0,\tilde{T}]$. Similarly, the lower semicontinuity holds for
$$\tilde I[w]:=\int_{0}^{\tilde{T}}\n w(t)-\ou(t)\en^{2}dt.$$
Passing to the limit as $k\to\infty$ in the cost functional, combining this with the assumed local continuity of $\varphi$, and using \eqref{J}, we conclude that the triple $(\tilde x (\cdot),\tilde u(\cdot),\tilde T)$ belongs to the prescribed
$W^{1,2}\times L^2 \times \R_+$-neighborhood of the given local minimizer
$(\ox(\cdot),\ou(\cdot),\oT)$ and satisfies the inequality
\begin{equation}\label{contr}
J[\tilde{x},\tilde{u},\tilde{T}]< J[\ox,\ou,\oT],
\end{equation}
which contradicts the fact that $(\ox(\oT),\ou(\oT),\oT)$ is a relaxed $W^{1,2}\times L^2\times\R_+$-local minimizer of $(P)$. Thus we arrive at \eqref{e:0.32} and complete the proof of the claimed convergence of optimal solutions from which the convergence of the
Lipschitz constants follows immediately.
\end{proof}\vspace*{-0.15in}

\section{Tools of Variational Analysis and Generalized Differentiation}\label{Var-Ana}\vspace*{-0.05in}

As mentioned above and will be clearly seen below, problems $(P)$ and $(P_k)$ unavoidably contain nonsmooth and nonconvex  components independently of smoothness and convexity of the given data. This is primary due to normal cone description of the sweeping differential inclusion and its discrete approximations generating {\em nonconvex graphical constraints}. To deal with such problems, we need to use suitable generalized differential constructions of variational analysis. It has been well recognized in the sweeping control theory that the limiting 
normal cone/coderivative/subdifferential notions for sets, set-valued mapping, and nonsmooth functions introduced by the second author are the most appropriate to derive necessary optimality conditions in sweeping optimal control. We also employ second-order subdifferentials of extended-real-valued functions that naturally appear in our derivation due to the very structure of the sweeping dynamics.\vspace*{-0.05in}

Let us start with recalling the employed {\em first-order} generalized differential constructions following the books \cite{m-book1,m-book2,rw}, where the reader can find proofs, further material, and bibliographies. The (Painlev\'e-Kuratowski) {\em outer limit} of a set-valued mapping/multifunction $F\colon\R^n\tto\R^m$ at $\ox$ with $F(\ox)\ne\emp$ is defined by
\begin{equation}\label{Pa}
\underset{x\to\ox}{\textrm{Lim sup }}F(x):=\big\{y\in\R^m\;\big|\;\exists\textrm{ sequences }\;x_k\to\ox,\;y_k\to y\;
\textrm{ such that }\;y_k\in F(x_k),\;k\in\N\big\}.
\end{equation}
Given a set $\O\subset\R^n$ and a point $\ox\in\Omega$, the (basic, limiting, Mordukhovich) {\em normal cone} to $\O$ at $\ox$ is
\begin{equation}\label{nor}
N(\ox;\O)=N_\O(\ox):=\underset{x\to\ox}{\textrm{Lim sup}}\;\hat N(x;\O),
\end{equation}
where the (Fr\'echet) {\em regular normal cone} to $\O$ at $x$ is defined by
\begin{equation*}
N(x;\O)=\hat N_\O(x):=\Big\{v\in\R^n\;\Big|\;\limsup_{u\st{\O}{\to}x}\frac{v,u-x}{\|u-x\|}\le 0\Big\}  
\end{equation*}
if $x\in\O$ and $\hat N(x;\O):=\emp$ otherwise, and where $u\st{\O}{\to}x$ means that $u\to x$ with $u\in\O$.\vspace*{-0.05in}

Associate with a set-valued mapping $F\colon\R^n\tto\R^m$ its domain and graph 
\begin{equation*}
\dom F:=\big\{x\in\R^n\;\big|\;F(x)\ne\emp\}\;\textrm{ and }\;\gph F:=\big\{(x,y)\in\R^n\times\R^m\;\big|\;y\in F(x)\big\}.
\end{equation*}
The {\em coderivative} of $F$ at $(\ox,\oy)\in\gph F$ is defined by
\begin{equation}\label{cod}
D^*F(\ox,\oy)(u):=\big\{v\in\R^n\;\big|\;(v,-u)\in N\big((\ox,\oy);\gph F\big)\big\}\;\textrm{ for }\;u\in\R^m.
\end{equation}
If $F\colon\R^n\to\R^m$ is single-valued and continuously differentiable $({\cal C}^1$-smooth) around $\ox$, then
\begin{equation*}
D^*F(\ox)(u)=\big\{\nabla F(\ox)^*u \big\}\;\textrm{ for all }\;u\in\R^m,
\end{equation*}
where $\nabla F(\ox)^*$ is the adjoint/transposed Jacobian matrix, and where $\oy=F(\ox)$ is omitted. 
When $F$ is single-valued and locally Lipschitzian around $\ox$, we have the scalarization formula:
\begin{equation}\label{scal}
D^*F(\ox)(u)=\partial\la u,F\ra(\ox)\;\mbox{ for all }\;u\in\R^m.
\end{equation}\vspace*{-0.05in}

Given an extended-real-valued function $\varphi\colon\R^n\to\oR:=(-\infty,\infty]$ with
\begin{equation}
\dom\varphi:=\big\{x\in\R^n\;\big|\;\vph(x)<\infty\big\}\;\textrm{ and }\;\epi\varphi:=\big\{(x,\alpha)\in\R^{n+1}\;\big|\;
\alpha\ge\varphi(x)\big\},
\end{equation}
the (first-order) {\em subdifferential} of $\ph$ at $\ox\in\dom\varphi$ is defined via the normal cone \eqref{nor} to the epigraph $\epi\ph$ by
\begin{equation}\label{sub1}
\partial\varphi(\ox):=\big\{v\in\R^m\;\big|\;(v,-1)\in N\big((\ox,\varphi(\ox));\epi\varphi\big)\big\}.
\end{equation}
Note that if $\ph(x):=\dd(x;\O)$ is the indicator function of $\O$ equal to $0$ for $x\in\O$ and $\infty$ otherwise, then $\partial\ph(\ox)=N(\ox;\O)$ for each $\ox\in\O$. Let us emphasize that the above normal cone, coderivative, and subdifferential constructions enjoy {\em full calculus} based on variational/extremal principles of variational analysis.\vspace*{-0.05in}

Following \cite{m92}, we define the {\em second-order subdifferential} (or {\em generalized Hessian}) of $\varphi\colon\R^n\to\oR$ at $\ox\in\dom\ph$ relative to $\oy\in\partial\varphi(\ox)$ via the coderivative of the first-order subdifferential mapping from \eqref{sub1} by
\begin{equation}\label{2nd}
\partial^2\varphi(\ox,\oy)(u):=(D^*\partial\varphi)(\ox,\oy)(u),\;u\in\R^n.
\end{equation}
If $\ph$ is ${\cal C}^2$-smooth around $\ox$ with the (symmetric) Hessian  $\nabla^2\varphi(\ox)$, then we have 
\begin{equation*}
\partial^2\varphi(\ox)(u)=\big\{\nabla^2\varphi(\ox)u\big\},\;u\in\R^n,
\end{equation*}
The sweeping process setting corresponds to the case where $\varphi(x):=\dd(x;\O)$. In this case,
\begin{equation}\label{2nd-sweep}
\partial^2\varphi(\ox,\ov)(u)=D^*N_\O(\ox,\ov)(u)\;\mbox{ for all }\;\ov\in N(\ox;\O)\;\mbox{ and }\;u\in\R^n.
\end{equation} 
The reader is referred to the book \cite{m24} for a comprehensive theory and various applications of \eqref{2nd} and related second-order constructions with extensive calculus rules and explicit evaluations/calculations of \eqref{2nd}, and particularly of \eqref{2nd-sweep}, via the given data of structured functions and sets that appear in broad frameworks of variational analysis, optimization, and control. In this paper, we employ the following second-order evaluations of the coderivative of the normal cone mappings $F^k_i$ from \eqref{F^k_i}, providing a brief proof with the references therein.\vspace*{-0.05in}

\begin{theorem}\label{Thr4}
Given $F^k_i$ in \eqref{F^k_i} with $C^k_i$ taken from \eqref{C_i}, let $M>0$ be sufficiently large. Suppose that $g$ is locally Lipschitzian around the points in question and that PLICQ \eqref{plicq} holds when $t=t^k_i$. Then for any $(x,u)\in C^k_i\times U$ and any
$w\in-g(x,u)+N(x;C^k_i)\cap M\B$, we have the upper estimate for the coderivative of $F^k_i$ in $(x,u)$ with omitting the dependence on $t\in[t^k_i,t^k_{i+1})$ for simplicity:
\begin{equation}\label{cod-est}
D^*F^k_i(x,u,w)(y)\subset\Big\{z\in \R^{n+d}\;\bigg|\;z\in\partial\la  y,-g\ra(x,u)+\Big(\underset{j\in I^i_{0k}(y),\;\alpha^j\in\R}\sum\alpha^j x^j_*(t^k_i)+\underset{j\in I^i_{>k}(y),\;\beta^j\geq 0}\sum\beta^j x^j_*(t^k_i),0\Big)\Big\},
\end{equation}
where $y\in\dom D^*N_{C^k_i}(x,w+g(x,u)-b\tau^k_i)$ as $b\in\B$, and where
\begin{equation}\label{I_k}
I^k_i(x):=\big\{j\;\big|\;\la x^j_{\ast}(t^k_i),x\ra \leq c^{ij}_k\big\},
\end{equation}
\begin{equation}\label{c56}
I^i_{0k}( y):=\big\{j\in I^k_i(x)\;\big|\;\la x^j_*(t^k_i),y\ra=0\big\}\;\textrm{ and }\;I^i_{>k}(y):=
\big\{j\in I^k_i(x)\;\big|\;\la x^j_*(t^k_i),y\ra>0\big\},\;y\in\R^n.
\end{equation} 
Furthermore, for each $t\in[t^k_i,t^k_{i+1})$ the equality in \eqref{cod-est} holds together with 
\begin{equation}\label{dom}
\begin{aligned}
\dom D^*F^k_i(x,u,w)=\Big\{y\;\Big|&\,\exists\,\lm^j\ge 0,\;b\in \B,\;j\in I^k_i(x)\;\mbox{ with }\;y+g(x,u)+b\tau^k_i=\underset{j\in I^k_i(x)}\sum\lm^jx^j_*(t^k_i),\\
&\mbox{and }\;\la x^j_{\ast}(t^k_i),y\ra =0\;\mbox{ if }\; 
\lm^j>0,\;\la x^j_{\ast}(t^k_i),y\ra \ge 0\;\mbox{ if }\;\lm^j=0 
\Big\}
\end{aligned}
\end{equation}
provided that the generating vectors $\{x^j_*(t^k_j)\;|\;j\in I(t^k_j,x)\}$ of the polyhedron $C^k_j$ are linearly independent.
\end{theorem}\vspace*{-0.05in}
\begin{proof}\; 
Employ the sum rule from \cite[Theorem~3.9(ii)]{m-book2}, where the sets
\begin{equation*}
S(x,u,w):=\big\{(a_1,a_2)\in\R^n\times\R^n\;\big|\;a_1\in F^k_{i1}(x,u),\;a_2\in F^k_{i2}(x,u),\;a_1+a_2=w\big\}
\end{equation*}
in the aforementioned theorem reduce in our case to the form
\begin{equation}\label{S}
S(x,u,w)=\Big\{\big(w+g(x,u),-g(x,u)\big)\in\R^n\times\R^n\;\Big|\;w\in -g(x,u)+N(x;C^k_i)-\tau^k_i\B\Big\}
\end{equation}
with $F^k_{i1}(x,u):= N(x;C^k_i)-\tau^k_i\B$ and $F^k_{i2}(x,u):=-g(x,u)$. The qualification condition in \cite[Theorem~3.9(ii)]{m-book2} is given in the coderivative form
\begin{equation*}
D^*F^k_{i1}(x,u,a_1)(0)\cap\big(-D^*F^k_{i2}(x,u,a_2)(0)\big)=\big\{(0,0)\big\},
\end{equation*}
and it holds by the assumed Lipschitz continuity of $g$ due to the necessity part of the coderivative criterion for Lipschitz continuity taken from \cite[Theorem~3.3]{m-book2}, which ensures that $D^*F^k_{i2}(x,u,a_2)(0)=\{(0,0)\}$. It follows from the coderivative sum rule in \cite[Theorem~3.9(ii)]{m-book2} applied to the sum $F^k_{i1}+F^k_{i2}$ that
\begin{eqnarray}\label{sum1}
D^*(F^k_{i1}+F^k_{i2})(x,u,w)(y)\subset\underset{(a_1,a_2)\in S(x,u,w)}{\bigcup}\Big(D^*F^k_{i1}(x,u,a_1)(y )+D^*F^k_{i2}(x,u,a_2)(y)\Big).
\end{eqnarray}
For each $k\in \N$, we use further yet another coderivative sum rule from \cite[Proposition~1.62]{m-book1} applied to the mapping $F^k_{i1}(x,u)$, which gives us $b\in\B$ such that  
\begin{equation*}
D^*F^k_{i1}(x,u,a_1)(y)=\big\{D^*N_{C^k_i}(x,w+g(x,u)-b\tau^k_i)(y)\;\big|\;y\in\dom D^*N_{C^k_i}\big(x,w+g(x,u)-b\tau^k_i\big)\big\}
\end{equation*}
with $\dom D^*F^k_{i1}(x,u,a_1)=\dom D^*N_{C^k_i}(x,w+g(x,u)-b\tau^k_i)$. Taking into account the scalarization formula \eqref{scal} for the Lipschitzian mapping $g(x,u)$ yields
$$
D^*(F^k_{i1}+F^k_{i2})(x,u,w)(y)\subset\big\{z\in\R^{n+d}\,\big|\,D^*N_{C^k_i}\big(x,w+g(x,u)-b\tau^k_i\big)(y)+\partial\la y,-g\ra(x,u)\big\}.
$$
Using the upper estimate of the coderivative of the normal cone mapping taken from \cite[Theorem~4.5]{hmn}, we get 
\begin{equation*}
\begin{array}{ll}
D^{\ast}N_{C^k_i}(x,w+g(x,u)-b\tau^k_i)(y)&\subset\Big(
\span\{x^j_*\;\big|\;j\in I^i_{0k}(y)\big\}+\cone\{x^j_*\;|\;j\in I^i_{>k}(y)\},0\Big)\\
&=\disp\Big(\underset{j\in I^i_{0k}(y),\alpha^j\in\R}\sum\alpha^j x^j_{\ast}(t^k_i)+\disp\underset{j\in I^i_{>k}(y),\beta^j\ge 0}\sum\beta^j x^j_{\ast}(t^k_i),0\Big),
\end{array}
\end{equation*}
which thus implies the inclusion in \eqref{cod-est}. Consider now the sets
$$
J:=\big\{j\in I(t^k_i, x)\;\big|\;\lambda^j>0\big\},\;\Upsilon(J):=\big\{j\in 
I(t^k_i, x)\;\big|\;\la x^{\ast}_j(t^k_i),x\ra =0\;\mbox{ for all }\;j\in C_J\big\}, 
$$
$$
\begin{aligned}
C_J:=\big\{x\in \R^n\;\big|\;\la x^{\ast}_j(t^k_i),x\ra =0\;\mbox{ for all }\;j\in J,\;\la x^{\ast}_j(t^k_i),x\ra \geq 0 \;\mbox{ for all }\;j\in \{1,\ldots,s\}\setminus J\big\}.	
\end{aligned}
$$
It follows from \cite[Theorem~4.5]{hmn} that 
$$
\Upsilon(J)\setminus J =\big\{j\in 
I(t^k_i, x)\;\big|\; \la x^{\ast}_j(t^k_i),x\ra =0\;\mbox{ for all }\;j\in C_J \;\mbox{ and }\;j\notin J\big\}. 
$$
This allows us to deduce from \cite[Theorem~4.5]{hmn} the domain representation
\begin{equation*}
\begin{aligned}
&\dom D^*N_{C^k_i}\big(x,u,w+g(x,u)-b\tau^k_i\big)=\Big\{y\;\Big|\;\exists\,\lm^j \geq 0,\;b\in\B,\;j\in I^k_i(x)\;\mbox{ with }\;w= y+g(x,u)+b\tau^k_i\\=&\underset{j\in I^k_i(x)}\sum\lm^j x^j_*(t^k_i),\;\mbox{ and }\;\la x^j_{\ast}(t^k_i),y\ra =0\;\mbox{ if }\;
\lm^j>0\;\mbox{ while}\;\la x^j_{\ast}(t^k_i),y\ra\geq 0\;\mbox{ if }\;\lm^j=0\Big\}
\end{aligned}
\end{equation*}
under the fulfillment of the imposed LICQ. This completes the proof of the theorem.
\end{proof}
\vspace*{-0.15in}

\section{Optimality Conditions for Discrete Approximations}\label{sec:NCO}\vspace*{-0.05in}

In this section, we derive necessary optimality conditions for problems $(P_k)$, $k\in\N$, 
as formulated in \eqref{d_a_p}--\eqref{F^k_i}. To accomplish this, we reduce each discrete-time problem $(P_k)$ to a nondynamic problem of nondifferentiable programming with increasingly many geometric and functional constraints. Employing necessary optimality conditions for the latter problem, established via the generalized differential constructions of Section~\ref{Var-Ana}, and then using generalized differential calculus rules and explicit second-order evaluations, we arrive at the desired optimality conditions for $(P_k)$ expressed via the given problem data. Here are the results.\vspace*{-0.07in}

\begin{theorem}\label{Thr5}
For each $k\in \N$, let $(\ox^k(\cdot),\ou^k(\cdot),\oT_k):=(\ox^k_0(\cdot),\ldots,\ox^k_{k}(\cdot),\ou^k_0(\cdot),\ldots,\ou^k_{k-1}(\cdot),\oT_k)$ be an optimal solution to problem $(P_k)$ under the general assumptions of Theorem~{\rm\ref{Thr4}}, and suppose that $\varphi$ is Lipschitz continuous around the point $(\ox^k(\oT_k),\oT_k)$. Then there exist a number $\mu^k_0\ge 0$ as well as vectors $\mu^k=(\mu^k_1,\ldots,\mu^k_s)\in\R^{s}_+$ and vectors $\{p^k_i\in\R^n\;|\;i=0,\ldots,k\}$ satisfying
the conditions
\begin{equation}\label{e5.1}
\mu^k_0+\n(\mu^k_1,\ldots,\mu^k_s)\en+\n \psi^k\en+\sum_{i=0}^{k}\n p^k_i\en\ne 0,
\end{equation}
\begin{equation}\label{e5.2}
\mu^k_i\big(\la x^{j}_\ast (t^k_k),x^k_k\ra-c_{j}(t^k_k)\big)=0,\quad j=1,\ldots,s,
\end{equation}
\begin{equation}\label{e5.4}
\begin{array}{ll}
&\disp\Big(\frac{p^k_{i+1}-p^k_i}{h^k},-\frac{\mu^k_0\xi^k_{iu}}{h^k},\frac{\mu^k_0\xi^k_{iy}}{h^k}-p^k_{i+1}\Big)\\
&\qquad\in\disp\Big(0,\frac{\psi^k_i}{h^k},0\Big)+N\Big(\Big(\ox^k_i,\ou^k_i,-\frac{\ox^k_{i+1}-\ox^k_i}{h^k}\Big);\gph F^k_i\Big),
\end{array}
\end{equation}
where $\psi^k=(\psi^k_0,\ldots,\psi^k_{k-1})$ with $\psi^k_i\in N(\ou^k_i;U)$ for $i=0,\ldots,k-1$, such that we have
\begin{equation}\label{e5.6}
\Big(-p^{k}_k-\disp\sum^{s}_{j=1}\mu^k_jx^j_*(t^k_k),\;\bar
{H}^{k}+2\mu^k_0(\oT-\oT_{k})+\mu^k_0\varrho_k\Big)\in\mu^k_0\partial\varphi(\bar{x}^{k}(
\oT_{k}),\oT_{k})+N\Big((\ox^k_k,\oT_k);\O^k_x\times\O^k_T\Big)
\end{equation}
with the data therein calculated by
\begin{equation}\label{t}
h^k=\dfrac{\oT_k}{k},\quad t^k_i=ih^k\;\mbox{ for }\; i=0,\ldots,k,
\end{equation}
\begin{equation}\label{Hk}
\bar{H}^k:=\dfrac{1}{k}\sum_{i=0}^{k-1}\la p^k_{i+1},y^k_i\ra,
\end{equation}
\begin{equation}\label{vrho}
\varrho_k:=\sum_{i=0}^{k-1}\left[\dfrac{i}{k}\n\(\frac{\ox^k_{i+1}-\ox^k_i}{h^k}-
\dot{\ox}(t^k_i),\ou^k_i-\ou(t^k_i)\)\en^2-\dfrac{i+1}{k}\n\(\frac{\ox^k_{i+1}-\ox^k_i}{h^k}-
\dot{\ox}(t^k_{i+1}),\ou^k_i-\ou(t^k_{i+1})\)\en^2\right]
\end{equation}
\begin{equation}\label{e5.5}
\xi^k_i:=\(\xi^k_{iu},\xi^k_{iy}\):=\left(\int_{t^k_i}^{t^k_{i+1}}
\(\ou^k_i-\ou(t)\) dt,\int_{t^k_i}^{t^k_{i+1}}\(\frac{\ox^k_{i+1}-\ox^k_i}{h^k}-\dot{\ox}(t)\)dt\right).
\end{equation}
\end{theorem}\vspace*{-0.05in}
\begin{proof}\;\;For any $k\in \N$, consider the nondynamic problem $(NP)$ with respect to the variable $$z:=(x^k_0,\ldots,x^k_k,u^k_0,\ldots,u^k_{k-1},y^k_0,\ldots,y^k_{k-1},\th),$$
where the starting point $x^k_0$ is fixed:
\begin{equation}\label{disc-cost}
\textrm{minimize }\;\phi_0(z):=\vph\big(x^k_k,\th\big)+(\th-\oT)^2+\sum_{i=0}^{k-1}\int_{\frac{i\th}{k}}^{\frac{(i+1)\th}{k}}\n\(y^k_i-
\dot{\ox}(t),u^k_i-\ou(t)\)\en^2dt
\end{equation}
subject to the functional and geometric constraints
\begin{equation*}
\omega_k(z):=\sum_{i=0}^{k-1}\int_{\frac{i\th}{k}}^{\frac{(i+1)\th}{k}}\n\(y^k_i-
\dot{\ox}(t),u^k_i-\ou(t)\)\en^2dt -\e\le 0,
\end{equation*}
\begin{equation*}
\omega_i(z):=\left\|\frac{x^k_{i+1}-x^k_i}{h^k}\right\| - L -1\le 0\;\text{ for all }\;i=0,\ldots, k-1,
\end{equation*}
\begin{equation*}
\phi_j(z):=\la x^{j}_\ast (t^k_k),x^k_k\ra-c_{j}(t^k_k)\le 0,\quad j=1,\ldots,s,
\end{equation*}
\begin{equation*}
g_i(z):=x^k_{i+1}-x^k_i-\dfrac{\th}{k} y^k_i=0,\quad i=0,\ldots,k-1,
\end{equation*}
\begin{equation*}
\Xi_i:=\big\{(x^k_0,\ldots,y^k_{k-1},\th)\;\big|\;-y^k_i\in F_i^k(x^k_i,u^k_i)\big\},\quad i=0,\ldots,k-1,
\end{equation*}
\begin{equation*}
z\in\Xi_k:=\big\{(x^k_0,\ldots,y^k_{k-1},\th)\;\big|\;x^k_0 \textrm{ is fixed},\;(x^k_k,\th)\in \O^k_x\times \O^k_T\},
\end{equation*}
\begin{equation*}
\Xi'_i:=\big\{(x^k_0,\ldots,y^k_{k-1},\th)\;\big|\; u^k_i\in U\big\},\quad i=0,\ldots,k-1.
\end{equation*}
It is clear that problems $(NP)$ and $(P_k)$ are equivalent. Denote by
$$
\oz:=\big(\ox^k_0,\ldots,\ox^k_k,\ou^k_0,\ldots,\ou^k_{k-1},\oy^k_0,\ldots,\oy^k_{k-1},\bar{\th}\big)
$$
the optimal solution to problem $(NP)$, where we drop the index $k$ if no confusion arises. Employing \cite[Theorem~6.5]{m-book2} gives us the necessary optimality conditions for $(NP)$ written as follows: there exist dual elements $\mu^k_0\ge 0$, $\mu^k=(\mu^k_1,\ldots,\mu^k_s)\in\R^{s}_+$, $\{p^k_i\in\R^n\;|\;i=1,\ldots,k\}$, and 
$$
z^*_i=\big(x^*_{0i},\ldots,x^*_{ki},u^*_{0i},\ldots,u^*_{(k-1)i},
y^*_{0i},\ldots,y^*_{(k-1)i},\th^*_i\big),\quad i=0,\ldots,k,
$$ 
not simultaneously equal to zero and such that we have
\begin{equation}\label{4.49}
z^*_i\in\left\{\begin{array}{ll}
N(\oz;\Xi_i\cap\Xi'_i)\;\textrm{ if }\;i\in\big\{0,\ldots,k-1\big\},\\
N(\oz;\Xi_k)\;\textrm{ if }\;i=k,
\end{array}\right.
\end{equation}
\begin{equation}\label{4.49'}
-z^*_0-\ldots-z^*_k\in\mu^k_0\partial\phi_0(\oz)+\sum_{j=1}^{s}\mu^k_j\nabla \phi_{j}(\oz)+
\sum_{i=0}^{k-1}\nabla g_i(\oz)^*p^k_{i+1},
\end{equation}
\begin{equation}\label{4.50}
\mu^k_j\phi_{j}(\oz)=0,\quad j=1,\ldots,s.
\end{equation}
Note that the inequality constraints in $(NP)$ defined by the functions $\omega_i$ as $i=0,\ldots, k$ are inactive for large $k$ due by Theorem~\eqref{Thr3}, and so the corresponding multipliers do not appear in the optimality conditions.\vspace*{-0.05in}

Next we claim that in inclusion \eqref{4.49}, the sets $\Xi_i$ and $\Xi'_i$ satisfy the qualification condition 
\begin{equation}\label{4.51}
N(\oz;\Xi_i)\cap\big(-N(\oz;\Xi'_i)\big)=\{0\},\quad i=0,\ldots,k-1.
\end{equation}
Indeed, fixing $z^*_i\in N(\oz;\Xi_i)\cap\big(-N(\oz;\Xi'_i)\big)$ and employing the limiting normal cone definition \eqref{nor} yield 
$$
z^*_i\in N(\oz;\Xi_i)= \underset{z\to\oz}{\textrm{Lim sup}}\;\hat{N}(z;\Xi_i),$$
which tells us that there are sequences $z_{im}\overset{\Xi_i}{\longrightarrow}\oz$, $v_{im}\overset{\Xi_i}{\longrightarrow}z_{in}$, and $z^*\to z^*_i$ satisfying
\begin{equation*}
\dfrac{\left \la z^*_{im}, v_{im}-z_{im} \right \ra}{\n v_{im}-z_{im}\en}\le 0\; \textrm{ for all }\; m \in \N.
\end{equation*}
This implies that $(\ox^k_i,\ou^k_i,-\oy^k_i\Big)\in \gph F^k_i$ with $-\oy^k_i=-\frac{\ox^k_{i+1}-\ox^k_i}{h^k}$ and $h^k= \frac{\th}{k}$, and that 
\begin{equation}\label{4.53}
\big(x^*_{ii},u^*_{ii},-y^*_{ii}\big)
\in N\Big(\Big(\ox^k_i,\ou^k_i,-\frac{\ox^k_{i+1}-\ox^k_i}{h^k}\Big);\gph F^k_i\Big),\quad i=0,\ldots,k-1,
\end{equation}
for the the corresponding components of $z^*_{i}$, while the other components of $z^*_{i}$ are equal to zero. Similarly, from $-z^*_i\in N(\oz;\Xi'_i)$ we get the inclusions
\begin{equation}\label{4.54}
-u^*_{ii}\in N(\ou^k_i;U),\quad i=0,\ldots,k-1. 
\end{equation}
Combining \eqref{4.53} and \eqref{4.54}, and then employing the same argument as above lead us to
$$
x^*_{ii}=0\;\mbox{ and }\;y^*_{ii}=0\;\mbox{ for }\;\ i=0,\ldots,k-1.
$$
Substitute the latter into \eqref{4.53} and using the coderivative definition \eqref{cod} give us
\begin{equation*}
(0,u^*_{ii})\in
 D^*F^k_i\Big(\ox^k_i,\ou^k_i,-\frac{\ox^{k+1}_i-\ox^k_i}{h^k}\Big)(0),\quad i=0,\ldots,k-1.
\end{equation*}
Next we employ the coderivative upper estimate \eqref{cod-est} in Theorem~\eqref{Thr4} valid under the uniform Slater condition at $y=0$. This yields $u^*_{ii}=0$ for all $i=0,\ldots, k-1$ and thus verifies the required qualification condition \eqref{4.51}.\vspace*{-0.05in}

Now we are in a position to apply the basic normal cone intersection rule from \cite[Theorem~2.16]{m-book2} to the sets $\Xi_i$ and $\Xi'_i$ in question, which tells us that
\begin{equation*}
z^*_i\in
N(\oz;\Xi_i)+N(\oz;\Xi'_i)\;\textrm{ if }\;i\in\big\{0,\ldots,k-1\big\}.
\end{equation*}
Hence there exist $\psi^k_i\in N(\ou^k_i;U)$ such that the first inclusion in \eqref{4.49} can be rewritten by
\begin{equation}\label{4.56}
\big(x^*_{ii},u^*_{ii}-\psi^k_i,-y^*_{ii}\big)
\in N\Big(\Big(\ox^k_i,\ou^k_i,-\frac{\ox^k_{i+1}-\ox^k_i}{h^k}\Big);\gph F^k_i\Big)\;\textrm{ for }
\;i=0,\ldots,k-1,
\end{equation}
while the other components of $z^*_i$ are equal to zero. For each $k\in\N$, the second inclusion in \eqref{4.49} reads as due 
\begin{equation*}
(x^{\ast}_{kk},\th^{\ast}_k)\in N((\ox^k_k,\bar{\th});\O^k_x\times \O^k_T).
\end{equation*}
Similarly, the only potentially nonzero component of $z^*_k$ is $x^{\ast}_{0k}$, which is determined by the normal cone to $\Xi_k$. Therefore, \eqref{4.49} tells us that
$$
\begin{aligned}
-z^*_0-\ldots-z^*_k=
\big(&-x^*_{00}-x^*_{0k},-x^*_{11},\ldots,-x^*_{k-1,k-1},-x^{\ast}_{kk},-u^*_{00},\ldots,-u^*_{k-1,k-1},\\
&-y^*_{00},-y^*_{11},\ldots,-y^*_{k-1,k-1},-\th^{\ast}_k\big).
\end{aligned}
$$\vspace*{-0.05in}

Next we calculate the sums on the right-hand side of \eqref{4.49'}. It follows from the constructions above that
\begin{eqnarray*}
\begin{aligned}
\Big(\sum_{j=1}^{s}\mu^k_j\nabla \phi_{j}(\oz)\Big)_{x^k_k}=&\sum_{j=1}^{s}\mu^k_j x_*^j(t^k_k),\\
\Big(\sum_{i=0}^{k-1}\nabla g_i(\oz)^*p^k_{i+1}\Big)_{x^k_i}=&
\left\{\begin{array}{lcl}-p^k_1\quad\mbox{if }\;i=0,\\
p^k_i-p^k_{i+1}\quad\mbox{if }\;i=1,\ldots,k-1,\\
p^k_k\quad\mbox{if }\;i=k,
\end{array}\right.\\
\Big(\sum_{i=0}^{k-1}\nabla g_i(\oz)^*p^k_{i+1}\Big)_{y^k_i}=&-\dfrac{\bar{\th}}{k}\( p^k_1,p^k_2,\ldots,p^k_k\),\\
\Big(\sum_{i=0}^{k-1}\nabla g_i(\oz)^*p^k_{i+1}\Big)_{\bar{\th}}=&-\dfrac{1}{k}\sum_{i=0}^{k-1}\la p^k_{i+1},y^k_i\ra.
\end{aligned}
\end{eqnarray*}
Applying the subdifferential sum rule to \eqref{disc-cost} gives us
\begin{eqnarray*}
\partial\phi_0(\oz)= \partial\vph(\ox^k_k,\bar{\th})+\nabla\rho(\oz)+2(0,\ldots,0,\bar{\th}-\bar{T)}\;\textrm{ with }\;
\rho(z):=\sum_{i=0}^{k-1}\int_{\frac{i\th}{k}}^{\frac{(i+1)\th}{k}}\n\(y^k_i-
\dot{\ox}(t),u^k_i-\ou(t)\)\en^2dt.
\end{eqnarray*}
Differentiating the function $\rho$ defined above, we easily have $\nabla_{x_i}\rho(\bar{z})=0$ for $i=0,\ldots,k$,
\begin{eqnarray*}
\nabla_{u_i}\rho(\bar{z})=2\int_{\frac{i\bar{\th}}{k}}^{\frac{(i+1)\bar{\th}}{k}}(\ou^k_i-\ou(t))dt,\;\nabla_{y_i}\rho(\bar{z})=2\int_{\frac{i\bar{\th}}{k}}^{\frac{(i+1)\bar{\th}}{k}}(\oy^k_i-\dot{\bar{x}}(t))dt
\end{eqnarray*}
for $i=0,\ldots,k-1$ together with the expression
\begin{eqnarray*}
\nabla_{\th}\rho(\bar{z})=\sum_{i=0}^{k-1}\left[\dfrac{i+1}{k}\n\(\oy^k_i-
\dot{\ox}\left(\dfrac{(i+1)\bar{\th}}{k}\right),\ou^k_i-\ou\left(\dfrac{(i+1)\bar{\th}}{k}\right)\)\en^2-\dfrac{i}{k}\n\(\oy^k_i-
\dot{\ox}\left(\dfrac{i\bar{\th}}{k}\right),\ou^k_i-\ou\left(\dfrac{i\bar{\th}}{k}\right)\)\en^2\right]
\end{eqnarray*}
for the $\th$-partial derivative of $\rho$ at $\bar z$ denoted by $\rho_{\th}$. It follows that 
$\mu^k_0\partial\phi_0(\oz)$ is the collection of elements 
$$
\mu^k_0(0,\ldots,0,\vth^k_k,\xi^k_{0u},\ldots,\xi^k_{(k-1)u},\xi^k_{0y},\ldots,\xi^k_{(k-1)y},\vth+\rho_{\th}+2(\bar{\th}-\oT)),
$$
where $(\vth^k_k,\vth)\in\partial\vph(\ox^k_k,\bar{\th})$, and where
$$(\xi^k_{iu},\xi^k_{iy}):=\left(2\int_{\frac{i\bar{\th}}{k}}^{\frac{(i+1)\bar{\th}}{k}}(\ou^k_i-\ou(t))dt,2\int_{\frac{i\bar{\th}}{k}}^{\frac{(i+1)\bar{\th}}{k}}(\oy^k_i-\dot{\bar{x}}(t))dt\right),\;\;i=0,\ldots,k-1.$$
Thus we split the inclusion in \eqref{4.49'} into the following equalities 
\begin{equation}\label{x_0}
-x^*_{00}-x^*_{0k}=-p^k_1,
\end{equation}
\begin{equation}\label{x_i}
-x^*_{ii}=p^k_i-p^k_{i+1},\quad i=1,\ldots,k-1,
\end{equation}
\begin{equation}\label{x_kk}
-x^{\ast}_{kk}=\mu^k_0\vth^k_k+p^k_k +\sum_{j=1}^{s}\mu^k_j x_*^j(t^k_k),
\end{equation}
\begin{equation}\label{u_i}
-u^*_{ii}=\mu^k_0\xi^k_{iu},\quad i=0,\ldots,k-1,
\end{equation}
\begin{equation}\label{y_i}
-y^*_{ii}=\mu^k_0\xi^k_{iy}-\dfrac{\bar{\th}}{k} p^k_{i+1},\quad i=0,\ldots,k-1,
\end{equation}
\begin{equation}\label{th_i}
-\th^{\ast}_k=\mu^k_0(\vth+\rho_{\th}+2(\bar{\th}-\oT))-\dfrac{1}{k}\sum_{i=0}^{k-1}\la p^k_{i+1},y^k_i\ra.
\end{equation}
Clearly, \eqref{4.50} yields \eqref{e5.2}. It follows from \eqref{x_i}, \eqref{u_i}, and \eqref{y_i} that
$$
\frac{x^*_{ii}}{h^k}=\frac{p^k_{i+1}-p^k_i}{h^k},\;\;\frac{u^*_{ii}}{h^k}=-\frac{1}{h^k}\mu^k_0\xi^k_{iu},\;
\textrm{ and }\;\frac{y^*_{ii}}{h^k}
=-\frac{1}{h^k}\mu^k\xi^k_{iy}+p^k_{i+1}\; \mbox{ for }\;i=0,\ldots , k-1,
$$
where $p^k_0:=x^{\ast}_{0k} $. Substituting this into \eqref{4.56}, we get \eqref{e5.4} and then deduce from \eqref{x_kk} and \eqref{th_i} that
$$
\dfrac{1}{k}\sum_{i=0}^{k-1}\la p^k_{i+1},y^k_i\ra-\mu^k_0\rho_{\th}+2\mu^k_0(\oT-\bar{\th})=\mu^k_0\vth+\th^{\ast}_k\;\text{ and }\;-p^k_k-\sum_{j=1}^{s}\mu^k_j x_*^j(t^k_k)=\mu^k_0\vth^k_k+x^{\ast}_{kk}.
$$
Since $(\vth^k_k,\vth)\in\partial\vph(\ox^k_k,\bar{\th})$ and 
$(x^{\ast}_{kk},\th^{\ast}_k)\in N((\ox^k_k,\bar{\th}),\O^k_x\times \O^k_T)$ , we arrive at
$$
\Big(-p^k_k-\sum_{j=1}^{s}\mu^k_j x_*^j(t^k_k),\dfrac{1}{k}\sum_{i=0}^{k-1}\la p^k_{i+1},y^k_i\ra-\mu^k_0\rho_{\th}+2\mu^k_0(\oT-\bar{\th}) \Big) \in\mu^k_0\partial\varphi(\ox^k_k,\bar{\th})+ N
\big((\ox^k_k,\bar{\th});\O^k_x\times\O^k_T\big)
$$
and thus justify the inclusion in \eqref{e5.6}. \vspace*{-0.05in}

It remains to verify the nontriviality condition \eqref{e5.1}. Indeed, supposing on the contrary that $\mu^k_i=0$ for $i=1,\ldots s$, $\psi^k_i=0$ for $i=0,\ldots, k-1$, and $p^k_i=0$ as $i=0,\ldots, k$ with taking into account that $p^k_0=x^*_{0k}$, we get $x^*_{0k}=p^k_0=0$. It follows from \eqref{x_0}, \eqref{x_i}, and \eqref{x_kk} that $x^*_{ii}=0$ for all $i=0,\ldots,k$. Using \eqref{u_i} tells us that $u^*_{ii}=0$ as $i=1,\ldots,k-1$. Furthermore, \eqref{y_i} yields $y^*_{ii}=0$ for all
$i=0,\ldots,k-1$. Remembering that all the components of $z^*_i$ different from $(x^*_{ii},u^*_{ii},y^*_{ii})$ are zeros for $i=0,\ldots,k-1$ ensures that $z^*_{i}=0$ for $i=0,\ldots,k-1$ and similarly $z^*_k=0$.
Therefore $z^*_i=0$ for all $i=0,\ldots,k$, which violates the nontriviality condition for $(NP)$ and hence completes the proof of theorem.
\end{proof}\vspace*{-0.03in}

The {\em discrete-time adjoint system} \eqref{e5.4} is actually expressed in terms of the {\em coderivative} \eqref{cod} of the  normal cone mapping $F^k_i$ from \eqref{F^k_i} that relates to the second-order subdifferential \eqref{2nd-sweep}. The next theorem employs the {\em second-order evaluation} of Theorem~\ref{Thr4} for such mappings to express the necessary optimality conditions for $(P_k)$ explicitly via the given data under our standing assumptions.\vspace*{-0.07in}

\begin{theorem}\label{Thr6} Let $(\ox^k(\cdot),\ou^k(\cdot),\oT_k)$ be an optimal solution to problem $(P_k)$ formulated in \eqref{d_a_p}--\eqref{re_6},
where the cost function $\ph$ is locally Lipschitzian around $(\ox^k(\bar{T_k}), \bar{T_k})$, and where the mappings $F^k_i$ are defined in \eqref{F^k_i}. Using the notation and assumptions of Theorem~{\rm\ref{Thr5}}, take $(\xi^k_{iu},\xi^k_{iy})$ from \eqref{e5.5}. Then for each $k\in\N$ sufficiently large, there exist dual elements $(\mu^k_0,\psi^k,p^k)$ as in Theorem~{\rm\ref{Thr5}} together with vectors $\eta^k_i\in\R^s_+$ for $i=0,\ldots,k$ and $\gg^k_i\in\R^s$ for $i=0,\ldots,k-1$ satisfying the nontriviality condition
\begin{equation}\label{ntc}
\mu^k_0+\n\eta^k_k\en+\sum_{i=0}^{k-1}\n p^k_i\en+\n \psi^k\en\ne 0,
\end{equation}
the primal-dual relationships given for all $i=0,\ldots,k-1$ and $j=1,\ldots,s$ by
\begin{equation}\label{diff}
b\tau^k_i-\frac{\ox^k_{i+1}-\ox^k_i}{h^k}+g(\ox^k_i,\ou^k_i)=\sum_{j\in I^k_i(\ox^k_i)}\eta^k_{ij} x^j_*(t^k_i),
\end{equation}
\begin{equation}\label{conx}
\begin{array}{ll}
&\left(\disp\frac{p^k_{i+1}-p^k_i}{h^k},-\frac{\mu^k_0\xi^k_{iu}}{h^k}-\frac{\psi^k_i}{h^k}\right)\in  \partial\Bigg\la-\dfrac{\mu^k_0\xi^k_{iy}}{h^k}+p^k_{i+1}, -g\Bigg\ra (\ox^k_i,\ou^k_i)\\
&+\(\disp\sum_{j\in I^i_{0k}(-\frac{\mu^k_0\xi^k_{iy}}{h^k}+p^k_{i+1})
\cup I^i_{>k}(-\frac{\mu^k_0\xi^k_{iy}}{h^k}+p^k_{i+1})}\gg^k_{ij} x^j_*(t^k_i),0\),
\end{array}
\end{equation}
and the transversality condition at the optimal endpoint and optimal ending time
\begin{equation}\label{p_kk}
(-p^{k}_k-\disp\sum^{s}_{j=1}\eta^k_{kj}x^j_*(t^k_k),\;\bar
{H}^{k}+2\mu^k_0(\oT-\oT_{k})+\mu^k_0\varrho_k)\in\partial\Big(\mu^k_0\varphi\Big)(\bar{x}^{k}(
\oT_{k}),\oT_{k})+N\((\ox^k_k,\oT_k);\O^k_x\times \O^k_T\),
\end{equation}
where $\psi^k_i\in N(\ou^k_i;U)$, $H_k$, and $\varrho_k$ as $i=0,\ldots,k-1$ are taken from Theorem~{\rm\ref{Thr5}}. We also have the complementary slackness condition \eqref{e5.2} together with 
\begin{equation}\label{eta}
\Big[\la x^j_\ast(t^k_i),\ox^k_i\ra<c_j(t^k_i)\Big]\Longrightarrow\eta^k_{ij}=0,
\end{equation}
\begin{equation}\label{93}
\left\{\begin{matrix}
\Big[j\in I^i_{>k}(-\frac{\mu^k_0\xi^k_{iy}}{h^k}+p^k_{i+1})\Big]\Longrightarrow\gg^k_{ij}\ge 0,\\
\Big[j\notin I^i_{0k}(-\frac{\mu^k_0\xi^k_{iy}}{h^k}+p^k_{i+1})\cup I^i_{>k}(-\frac{\mu^k_0\xi^k_{iy}}{h^k}+p^k_{i+1})\Big]
\Longrightarrow\gg^k_{ij}=0.
\end{matrix}\right.
\end{equation}
\begin{equation}\label{94}
\[\la x^j_\ast(t^k_i),\ox^k_i\ra<c_j(t^k_i)\]\Longrightarrow\gg^k_{ij}=0,
\end{equation}
\begin{equation}\label{eta1}
\big[\la x^j_\ast(t^k_i),\ox^k_k\ra < c_j(t^k_i)\big]\Longrightarrow\eta^k_{kj}=0
\end{equation}
for $i=0,\ldots,k-1$ and $j=1,\ldots,s$. If the vectors $\{x^j_*(t^k_i)|\;j\in I^k_i(\ox^k_i)\}$ are linearly independent, then we get
\begin{equation}\label{96}
\eta^k_{ij}>0\Longrightarrow\Big[\Big\la x^j_\ast(t^k_i),-\frac{\mu^k_0\xi^k_{iy}}{h^k}+p^k_{i+1}\Big\ra=0\Big]
\end{equation}
along with the enhanced nontriviality condition
\begin{equation}\label{entc}
\mu^k_0+\|\eta^k_k\|+\|p^k_0\|+\|\psi^k\|\ne 0.
\end{equation}
\end{theorem}\vspace*{-0.05in}
\begin{proof}\;\;
By the coderivative definition \eqref{cod} and the necessary condition \eqref{e5.4} in Theorem~{\rm\ref{Thr5}}, we get 
\begin{equation}\label{cod-disc}
\Big(\frac{p^k_{i+1}-p^k_i}{h^k},-\frac{\mu^k_0\xi^k_{iu}}{h^k}-\frac{\psi^k_i}{h^k}\Big)\in
 D^*F^k_i\Big(\ox^k_i,\ou^k_i,-\frac{\ox^k_{i+1}-\ox^k_i}{h^k}\Big)\(-\frac{\mu^k_0\xi^k_{iy}}{h^k}+p^k_{i+1}\)
\end{equation}
for all $i=0,\ldots,k-1$. It follows from the inclusion
\begin{equation}
\tau^k_i\B-\frac{\ox^k_{i+1}-\ox^k_i}{h^k}+g(\ox^k_i,\ou^k_i)\subset N(\ox^k_i;C^k_i)\;\textrm{ for }\;i=0,\ldots,k-1
\end{equation}
that there exist vectors $\eta^k_i \in\R^s_+$ as $i=0,\ldots,k-1$ and $b\in\B$ such that
\begin{equation*}
b\tau^k_i-\frac{\ox^k_{i+1}-\ox^k_i}{h^k}+g(\ox^k_i,\ou^k_i)=\sum_{j\in I^k_i(\ox^k_i)}\eta^k_{ij} x^j_*(t^k_i),
\end{equation*}
which verifies \eqref{diff} and \eqref{eta}. Applying the upper coderivative estimate \eqref{cod-est} in Theorem~{\rm\ref{Thr4}} at $x:=\ox^k_i$, $u:=\ou^k_i$,
$w:=-\frac{\ox^k_{i+1}-\ox^k_i}{h^k}$, and $y:=-\frac{\mu^k_0\xi^k_{iy}}{h^k}+p^k_{i+1}$ for $i=0,\ldots,k-1$ gives us 
$\gg^k_i\in\R^s$ and the relationships
$$
\Big(\frac{p^k_{i+1}-p^k_i}{h^k},-\frac{\mu^k_0\xi^k_{iu}}{h^k}-\frac{\psi^k_i}{h^k}\Big)
\in \partial\Bigg\la -\frac{\mu^k_0\xi^k_{iy}}{h^k}+p^k_{i+1},-g\Bigg\ra (\ox^k_i,\ou^k_i)+\(\disp\sum_{j\in I^i_{0k}(-\frac{\mu^k_0\xi^k_{iy}}{h^k}+p^k_{i+1})
\cup I^i_{>k}(-\frac{\mu^k_0\xi^k_{iy}}{h^k}+p^k_{i+1})}\gg^k_{ij} x^j_*(t^k_i),0\),$$
$$\psi^k_i\in N(\ou^k_i;U)\;\textrm{ as }\;i=0,\ldots,k-1.
$$
This justifies the conditions in \eqref{conx}, \eqref{93}, and \eqref{94}. Assuming now that the generating vectors $\{x^j_*(t^k_i)\;|\;j\in I^k_i(\ox^k_i)\}$ are linearly independent, we arrive at \eqref{diff} by applying the domain calculation for $D^*F^k_i(\ox^k_i,-\frac{\ox^k_{i+1}-\ox^k_i}{h^k}+g(\ox^k_i,\ou^k_i))$ in Theorem~\ref{Thr4}. This yields the implications
\begin{equation*}
\eta^k_{ij}>0\Longrightarrow\Big[\Big\la x^j_*(t^k_i),-\frac{\mu^k_0\xi^k_{iy}}{h^k}+p^k_{i+1}\Big\ra=0\Big]
\end{equation*}
and justifies \eqref{96}. Denote $\eta^k_k:=(\mu^k_1,\ldots,\mu^k_s)$ with $\mu^k_i$ from Theorem~\ref{Thr5} and hence get $\eta^k_i\in\R^s_+$ as $i=0,\ldots,k$. Implications \eqref{eta1} follow directly from \eqref{e5.2} and the definition of
$\eta^k_k$. It remains to verify the enhanced nontriviality \eqref{entc}. Supposing the contrary that $\mu^k_0=0$, $\eta^k_k=0$, $p^k_0=0$, and $\psi^k=0$ yields $p^k_{i+1}=0$ as $i=0,\ldots,k-1$ by \eqref{conx}. This shows that the nontriviality condition \eqref{ntc} fails and so completes the proof. 
\end{proof}\vspace*{-0.15in}

\section{Necessary Conditions in Sweeping Optimal Control}\label{nc-sweep}
\setcounter{equation}{0}\vspace*{-0.05in}

In this section, by passing to the limit as $k\to\infty$ from the necessary optimality conditions of Theorem~\ref{Thr6} and combining this with Theorem~\ref{Thr3} and the tools of generalized differentiation discussed in Section~\ref{Var-Ana}, we arrive at our main results providing necessary optimality conditions for the designated relaxed $W^{1,2}\times L^2\times\R_+$--local minimizer in $(P)$ under the standing assumptions imposed in Section~\ref{intro}.\vspace*{-0.05in}

Recall that Motzkin's theorem of the alternative ensures the normal cone representation
\begin{equation}\label{motzkin}
N\(x(t);C(t)\)=\Big\{\sum_{j\in I(t,x(t))}\al^j(t)x^j_*(t)\;\Big|\;\al^j(t)\ge 0\Big\}.
\end{equation}
Having this in hands, we introduce the following notion important for our subsequent considerations. \vspace*{-0.1in}

\begin{definition}\label{active cone}
Let $x(\cdot)$ be a solution to the controlled sweeping process \eqref{Problem}, i.e.,
\begin{equation*}
-\dot{x}(t)=\sum_{j=1}^s\eta^j(t)x^j_*(t)-g\big(x(t),u(t)\big)\;\textrm{ for all }\;t\in [0,T),
\end{equation*}
where $\eta_j\in L^2([0,T];\mathbb{R}^+)$ and $\eta_j(t) =0$ for a.e. $t$ such that $j\not\in I(t,x(t))$. We say that the normal cone to
$C(t)$ is {\sc active along} $x(\cdot)$ on the set $E\subset[0,T]$ if for a.e. $t\in E$ and all $j\in I(t,x(t))$ it holds that
$\eta_j(t)>0$. Denoting by $E_0$ the largest subset of $[0,T]$ where the normal cone to $C(t)$ is active along $x(\cdot)$\footnote{ Here we mean the union of the density points of $E$ such that the normal cone to $C(t)$ is active along $E$.}, we say that it is simply active along $x(\cdot)$ provided that $E_0=[0,T]$.
\end{definition}\vspace*{-0.1in}

Note that the requirement that the normal cone is active along a trajectory of \eqref{Problem} distinguishes \eqref{Problem} from trajectories of the classical controlled ODE $\dot{x}=g(x,u)$ which satisfy the {\em tangency condition} $\langle g(x(t),u(t)), x_\ast^j(t)\rangle \le 0$ for a.e. $t$ such that $j\in I(t,x(t))$ and thus $x(t)\in C(t)$ for all $t\in E$. Observe also that the above requirement holds automatically when the set $I(t,x(t))$ is empty, i.e., when $x(t)$ stays in the interior of $C(t)$.\vspace*{-0.1in}

\begin{theorem}\label{Thr7}
Let $(\ox(t),\ou(t),\oT),\;0\leq t\leq\oT$, be a relaxed $W^{1,2}\times L^2\times\R_+$--local minimizer to problem $(P)$. In addition to {\rm(H1)}, {\rm(H3)}, and {\rm(H4)}, suppose that LICQ holds along $\ox(t)$, $t\in[0,T]$, and that $\ph$ is locally Lipschitzian around $(\bar{x}(\oT),\oT)$. Then there exist a multiplier $\mu\ge 0$, a nonnegative vector measure 
$\gg_>=(\gg_>^1,\ldots,\gg_>^s)\in C^*([0,\oT];\R^s)$, and a signed vector measure 
$\gg_0=(\gg_0^1,\ldots,\gg_0^s)\in C^*([0,\oT];\R^s)$ together with adjoint arcs $p(\cdot)\in W^{1,2}([0,\oT];\R^n)$ and $q(\cdot)\in BV([0,\oT];\R^n)$ such that the following conditions hold:\\[1ex]
$\bullet$ The {\sc sweeping trajectory representation}
\begin{equation}\label{37}
-\dot{\ox}(t)=\sum_{j=1}^s\eta^j(t)x^j_*(t)-g\big(\ox(t),\ou(t)\big)\;\textrm{ for a.e. }\;t\in [0,\oT),
\end{equation}
where the functions $\eta^j(\cdot)\in L^2([0,\oT]);\R_+)$ are uniquely determined for a.e.\ $t\in[0,\oT)$ by representation \eqref{37}.\\[1ex]
$\bullet$ The {\sc adjoint arc inclusion}
\begin{equation}\label{c:6.6}
(-\dot{p}(t), \psi(t))\in \co\partial\la q(t), g\ra (\ox(t),\ou(t))\;\textrm{ for a.e. }\;t\in[0,\oT],
\end{equation}
where the subdifferential is taken with respect to $(x,u)$, where $\psi(\cdot)\in L^2([0,\oT];\R^d)$ satisfies the 
\begin{equation}\label{co}
\psi (t)\in \co N(\ou(t);U)\;\textrm{ for a.e. }\;t\in[0,\oT],
\end{equation}
where the right continuous representative of $q(\cdot)$ is
given by
\begin{equation}\label{c:6.9}
q(t)=p(t)-\int_{(t,\oT]} \sum_{j=1}^sd\gg^j(\tau)x^j_*(\tau)
\end{equation}
for a.e. $t\in[0,\oT]$ except at most a countable subset, and where $\gg:=\gg_>+\gg_0$. Moreover, $p(\oT)=q(\oT)$.\\[1ex]
$\bullet$ The {\sc tangential maximization condition}: if the normal cone \eqref{nor} is generated as
\begin{equation}\label{dua}
N(\ou(t);U)=T^*(\ou(t);U):=\big\{v\in\R^n\big|\;\la v,u\ra\le 0\;\mbox{ for all }\;u\in T\big(\ou(t);U\big)\big\}
\end{equation}
by some tangent set $T(\ou(t);U)$ associated with $U$ at $\ou(t)$, then we have
\begin{equation}\label{c:6.6'}
\big\la \psi(t),\ou(t)\big\ra=\max_{u\in T(\ou(t);U)}\big\la \psi(t),u\big\ra\;\textrm{ for a.e. }\;t\in[0,T].
\end{equation}
In particular, the {\sc global maximization condition}
\begin{equation}\label{max}
\big\la \psi(t),\ou(t)\big\ra=\max_{u\in U}\big\la \psi(t),u\big\ra\;\textrm{ for a.e. }\;t\in[0,\oT]
\end{equation}
is satisfied provided that the control set $U$ is convex.\\[1ex]
$\bullet$ The {\sc dynamic complementary slackness conditions}
\begin{equation}\label{41}
\big\la x^{j}_\ast(t),\ox(t)\big\ra<c_j(t)\Longrightarrow\eta^j(t)=0\;\mbox{ and }\;\eta^j(t)>0\Longrightarrow\;\big\la x^{j}_\ast(t),q(t)\big\ra=0
\end{equation}
for a.e.\ $t\in[0,\oT]$ and all indices $j=1,\ldots,s$ if in addition LICQ at $\ox(t)$ is assumed.\\[1ex]
$\bullet$ The {\sc transversality conditions at the optimal final time}: there exist numbers $\eta^j(\oT)\ge 0$ whenever $j\in I(\oT,\ox(\oT))$ ensuring the relationships
\begin{equation}\label{42}
(-p(\oT)-\sum_{j\in I(\oT,\ox(\oT))}\eta^j(\oT) x^{j}_\ast(\oT),\bar{H})\in\mu \partial\varphi(\bar{x}(\oT),\oT)+N((\ox(\oT),\oT); \O_x\times \O_T),\;\;\eta^j(\oT)>0\Longrightarrow j\in I(\oT,\ox(\oT)),
\end{equation}
where $\bar{H}:={\oT}^{-1}\int_{0}^{\oT}\langle p(t),\dot{\bar{x}}(t)\rangle dt$ is a characteristic of the optimal time.\\[1ex]
$\bullet$ The {\sc endpoint complementary slackness conditions}
\begin{equation}\label{41a}
\big\la x^{j}_\ast(\oT),\ox(\oT)\big\ra<c_j(\oT)\Longrightarrow\eta^j(\oT)=0\;\mbox{ for all }\;j\in I(\oT,\ox(\oT))
\end{equation}
with the nonnegative numbers $\eta^j(\oT)$ taken from \eqref{42}.\\[1ex]
$\bullet$ The {\sc nonatomicity condition:} If $t\in[0,\oT)$ and $\la x^{j}_\ast(t),\ox(t)\ra<c_j$ for all $j=1,\ldots,s$, then there
exists a neighborhood $V_t$ of $t$ in $[0,\oT)$ such that $\gg^j_0(V)=\gg^j_>(V)=0$ for all Borel subsets $V$ of $V_t$. In particular, $\mathrm{supp}(\gamma^j_>)$ and $\mathrm{supp}(\gamma^j_0)$ are contained in the set $\{t\;|\;j\in I(t,\ox(t))\}$.\\[1ex]
$\bullet$ The {\sc general nontriviality condition}
\begin{equation}\label{e:83}
(\mu,p,\|\gamma_0\|_{TV},\|\gg_>\|_{TV})\ne 0,
\end{equation}
accompanied by the {\sc support condition}
\begin{equation}\label{suppcond}
\mathrm{supp}(\gamma_>)\cap \mathrm{int}(E_0)=\emptyset,
\end{equation}
that holds provided the normal cone is active on a set with nonempty interior.\\[1ex]
$\bullet$ The {\sc enhanced nontriviality}: we have $\mu=1$ for the cost function multiplier in \eqref{42} provided that $\la x^j_\ast(t),\ox(t)\ra<c_j(t)$ for all $t\in[0,\oT]$ and all indices $j=1,\ldots,s$.
\end{theorem}\vspace*{-0.1in}
{\bf Proof.} Given the local minimizer $(\ox(\cdot),\ou(\cdot),\oT)$ for $(P)$, construct the discrete-time problems $(P_k)$ for which
the existence of the optimal solutions $(\ox_k(\cdot),\ou_k(\cdot),\oT_k)$ is derived in Proposition~\ref{Thr1} and the convergence to $(\ox(\cdot),\ou(\cdot),\oT)$ is obtained by Theorem~\ref{Thr3}. We deduce each of the claimed necessary conditions in $(P)$ by passing to the limit from those in Theorem~\ref{Thr6}. Let us split the derivation into the {\em five steps} as follows.\\[1ex]
{\bf Step~1:} {\em Sweeping arc representation and dynamic complementary slackness.} Let us begin with verifying \eqref{37} together with the first implication in \eqref{41}. Using \eqref{e5.5}, define the functions $\xi^k(t)=(\xi^k_y(t),\xi^k_u(t))$ on $[0,\oT_k]$ by
$$
\xi^k_y(t):=\frac{\xi^{k}_{iy}}{h^k}\;\mbox{ and }\;\xi^k_u(t):= \frac{\xi^{k}_{iu}}{h^k}\;\textrm{ for }\;t\in[t^k_i,t^{k}_{i+1})\textrm{ and }\;i=0,\ldots,k-1.
$$
We can easily check by the constructions that
\begin{eqnarray*}
\int_0^{\oT_k}\n\xi^k_y(t)\en^2dt&=&\sum_{i=0}^{k-1}\frac{\Big\|\xi^k_{iy}\Big\|^2}{h^k}\le\frac{1}{h^k}\sum_{i=0}^{k-1}
\Big(\int_{t^k_i}^{t^k_{i+1}}\n\dot{\ox}(t)-\frac{\ox^k_{i+1}-\ox^k_i}{h^k}\en dt\Big)^2\\
&\le&\sum_{i=0}^{k-1}\int_{t^k_i}^{t^k_{i+1}}\Big\|\dot{\ox}(t)-\frac{\ox^k_{i+1}-\ox^k_i}{h^k}\Big\|^2dt=
\int_0^{\bar{T_k}}\n\dot{\ox}(t)-\dot{\ox}^k(t)\en^2dt.
\end{eqnarray*}
If follows from the strong convergence $(\ox^k(\cdot),\ou^k(\cdot))\to(\ox(\cdot),\ou(\cdot))$ in Theorem~\ref{Thr3} that
\begin{equation}\label{c:6.14}
\int_0^{\oT_k}\n\xi^k_y(t)\en^2dt\le\int_0^{\oT_k}\n\dot{\ox}(t)-\dot{\ox}^k(t)\en^2dt\to0\;\textrm{ as }\;k\to\infty.
\end{equation}
This shows that a subsequence of $\{\xi^k_y(t)\}$ (no further relabeling) converges to zero a.e.\ on $[0,\oT]$. We similarly get
\begin{eqnarray*}
\int_0^{\oT_k}\Big\|\xi^k_u(t)\Big\|^2dt&=&\sum_{i=0}^{k-1}\frac{\Big\|\xi^k_{iu}\Big\|^2}{h^k}\le\frac{1}{h^k}
\sum_{i=0}^{k-1}
\Big(\int_{t^k_i}^{t^k_{i+1}}\n\ou^k_i-\ou(t)\en dt\Big)^2\\
&\le&\sum_{i=0}^{k-1}\int_{t^k_i}^{t^k_{i+1}}\n\ou^k_i-\ou(t)\en^2dt=\int_0^{\oT_k}\n\ou^k(t)-\ou(t)\en^2dt,
\end{eqnarray*}
which being combined with Theorem~\ref{Thr3} ensures that
\begin{equation}\label{c:6.14'}
\int_0^{\oT_k}\n\xi^k_u(t)\en^2dt\le\int_0^{\oT_k}\n\ou^k(t)-\ou(t)\en^2dt\to 0\;\textrm{ as }\;k\to\infty,
\end{equation}
and hence $\xi^k_u(t)\to 0$ for a.e.\ $t\in[0,\oT]$ as $k \to \infty$ along a subsequence. By the choice of $x^{j}_\ast(t^k_i)$ and the convergence in Theorem~\ref{Thr3}, we deduce from the robustness of LICQ along $\ox(\cdot)$ and that the vectors $\{x^{j}_\ast(t^k_i)\;|\;j\in I^k_i(\ox^k_i)\}$ are linearly independent for each $i=1,\ldots,k$ and all $k\in\N$ sufficiently large.\vspace*{-0.05in}

Taking $\eta^k_{i}\in\R^s_+$ from Theorem~\ref{Thr6}, construct the piecewise constant functions $\eta^k(\cdot)$ on
$[0,\oT_k]$ by $\eta^k(t):=\eta^k_{i}$ for $t\in[t^k_i,t^k_{i+1})\;\textrm{ as }\;i=0,\ldots,k-1$, where $\eta ^k_i\in \R^s_+$. It follows from \eqref{diff} that 
\begin{equation}\label{c:51}
-\dot{\ox}^k(t)=\sum_{j=1}^s\eta^k_j(t)x^j_*(t)-g\big(\ox^k(t),\ou^k(t)\big)- \tau^k(t)\B
\;\textrm{ whenever }\;
t\in[t^k_i,t^k_{i+1}),\quad k\in\N.
\end{equation}
By passing to the limit in \eqref{c:51} and employing the normal cone representation \eqref{motzkin} give us $-\dot{\ox}(t)\in N(\ox(t);C(t))-g(\ox(t),\ou(t))$ for a.e.\ $t\in[0,\oT]$, where the normal cone mapping $t\mapsto N(t,C(t))$ is measurable on $[0,\oT]$ by \cite[Theorem~14.26]{rw}. It follows for the measurable selection result of \cite[Corollary~14.6]{rw} that there is measurable 
\begin{equation*}
(\al^j)_{j\in I(t,\ox(t))}(t)\mapsto\dot{\ox}(t)-g\big(\ox(t),\ou(t)\big)+\Big\{\sum_{j\in I(t,\ox(t))}\al^j(t) x^j_*(t)\;\Big|\;\al^j(t)\ge 0\Big\},
\end{equation*}
which is a.e.\ uniquely determined due to the imposed LICQ along $\ox(t)$.
Denote
\begin{equation}\label{eta-mes}
\eta^j(t):=\left\{\begin{array}{ll}
\al^j(t)&\mbox{for }\;j\in I(t,\ox),\\
0&\mbox{otherwise}
\end{array}\right.
\end{equation}
on $[0,\oT)$ and observe that each $\eta^j(t)$ as $j\in {1, \ldots, s}$ belongs to $L^2([0,\oT];\R_+)$ by combining  \eqref{eta-mes}, $\dot \ox(\cdot)\in L^2([0,\oT];\R^n)$, and (H3). This verifies the claimed sweeping arc representation \eqref{37} for a.e. $t\in [0,\oT)$. To define $\eta^j(t)$ at the endpoint $t=\oT$, take $\eta^k(\oT):=\eta^{k}_k$ from the optimality conditions for discrete approximations in Theorem~\ref{Thr6} and deduce from the nontriviality conditions in \eqref{ntc} after their normalization that the sequence $\{\eta^k_k\}$ converges, along a subsequence, to some vector $(\eta^{1}(\oT),\ldots,\eta^{k}(\oT))$. Combining \eqref{c:51} and \eqref{37} ensures that
\begin{equation*}
\dot{\ox}(t)-\dot{\ox}^k(t)=\sum_{j=1}^s\big[\eta^k_j(t)-\eta_j(t)\big]x^j_*(t)+g\big(\ox(t),\ou(t)\big)-g\big(\ox^k(t),
\ou^k(t)\big)- \tau^k(t)\B
\end{equation*}
for $t\in[t^k_i,t^k_{i+1})$ and $i=0,\ldots,k-1$, $k\in \N$.
It can be easily seen from the above that
\begin{equation*}
\Big\|\sum_{j=1}^s\big[\eta_j(t)-\eta^k_j(t)\big]x^j_*(t)\Big\|_{L^2}\le\n\dot{\ox}^k(t)-\dot{\ox}(t)\en_{L^2}+\n g\big(\ox(t),\ou(t)\big)-g\big(\ox^k(t),\ou^k(t)\big)\en_{L^2}+\tau^k(t)\B
\end{equation*}
whenever $t\in(t^k_i,t^k_{i+1})$. Letting $k\to \infty$ and combining this with Theorem~\ref{Thr3} give us
\begin{equation*}
\sum_{j\in I(t,\ox)}\big[\eta_j(t)-\eta^k_j(t)\big]x^j_*(t)\to 0\;\textrm{ for a.e. }\;t\in[0,\oT]\;\textrm{ as } k\to\infty,
\end{equation*}
and leads therefore to the a.e.\ convergence $\eta^k(t)\to\eta(t)$ on $[0,\oT]$ by the assumed LICQ. Then verifies by \eqref{eta} the first complementarity implication in \eqref{41}.\\[1ex]
{\bf Step~2:} {\em Adjoint arcs and maximization conditions.} First we construct extensions of the discrete adjoint variables in Theorem~\ref{Thr6} to the entire interval $[0,\oT_k]$ for each $k$. Define $q^k(t)$ by
extending $p^k_i$ piecewise linearly on $[0,\oT_k]$ with $q^k(t^k_i):=p^k_i$ for $i=0,\ldots,k$.
Construct further $\gg^k(t)$ and $\psi^k(t)$ on $[0,\oT_k]$ by
\begin{equation}\label{c:6.25}
\gg^k(t):=\gg^k_i,\quad\psi^k(t):=\frac{1}{h^k}\psi^k_i\;\textrm{ for }\;t\in[t^k_i,t^k_{i+1})\;\textrm{ and }\;i=0,
\ldots,k-1
\end{equation}
with $\gg^k(\oT_k):=0$ and $\psi^k(\oT_k):=0$. Consider the functions
\begin{equation*}
\nu^k(t):=\max\big\{t^k_i\;\big|\;t^k_i\le t,\;0\le i\le k-1\big\}\;\textrm{ for all }\;t\in[0,\oT_k],\quad k\in\N,
\end{equation*}
and then obtain from $(\ref{conx})$ that, respectively,
\begin{equation}\label{conx'}
\begin{array}{ll}
&(\dot{q}^k(t), -\mu^k_0\xi^k_u(t)-\psi^k(t))\in
\partial\la -\mu^k_0\xi^k_y(t)+q^k(\nu^k(t)+h^k), -g \ra  \big(\ox^k(\nu^k(t)),\ou^k(\nu^k(t))\big)\\
& + \bigg(\disp\sum_{j\in I_{0}(\nu^k(t),-\mu^k_0\xi^k_y(t)+q^k(\nu^k(t)+h^k ))\cup I_{>}(\nu^k(t),-\mu^k_0\xi^k_y(t)+q^k(\nu^k(t)+h^k))}\gg^k_{j}(t) x^j_*(\nu^k(t)), 0\bigg),
\end{array}
\end{equation}
where $I_0(t,y):= \{j\in I(t,x)\;|\;\la x^j_{\ast}(t),y \ra =0\}$, and $I_>(t,y):= \{j\in I(t,x)\;|\;\la x^j_{\ast}(t),y\ra >0\}$
for every $t\in(t^k_i,t^k_{i+1})$ and $i=0,\ldots,k-1$. Define the adjoint arcs $p^k(\cdot)$ on $[0,\oT_k]$ by
\begin{equation}\label{c:6.29}
p^k(t):=q^k(t)+\int_t^{\oT_k}\sum_{j=1}^s\gg^k_j(\tau)x^j_\ast(\tau)d\tau\;\textrm{ for each }\;t\in[0,\oT_k].
\end{equation}
This shows that $p^k(\oT_k)=q^k(\oT_k)$, and furthermore
\begin{equation}\label{c:6.30}
\dot{p}^k(t)=\dot{q}^k(t)-\sum_{j=1}^s\gg^k_j(t)x^j_\ast(t)\;\textrm{ for a.e. }\;t\in[0,\oT_k].
\end{equation}
For simplicity of presentation, suppose in what follows that $g(x,u)$ is continuously differentiable in $x$} around the given local minimizer $\ox(t)$ uniformly in $u\in U$ and $t\in[0,\oT]$. This allows us to avoid technicalities in the subsequent exposition without increasing the length of the paper, while the reader can check that the arguments below hold under the general Lipschitzian assumption in (H3) due to the uniform boundedness and robustness of the basic subdifferential in \eqref{conx'} for locally Lipschitzian functions; see \cite{m-book2}. Taking this into account, we deduce from \eqref{conx'}, \eqref{93}, \eqref{c:6.30}, and the index definitions in \eqref{c56} that
\begin{equation}\label{c:59}
\dot{p}^k(t)=-\nabla_x g\big(\ox^k(\nu^k(t)),\ou^k(\nu^k(t))\big)^*\big(-\mu^k_0\xi^k_y(t)+q^k(\nu^k(t)+h^k)\big)
\end{equation}
for all $t\in(t^k_i,t^k_{i+1})$ and $i=0,\ldots,k-1$. Next we define the vector measures $\gg^k$ on $[0,\oT_k]$ by
\begin{equation}\label{c:6.34}
\underset{B_k}{\int}d\gg^k:=\underset{B_k}{\int}\sum_{i=0}^{k-1}\gg^k(t)\mathbbm{1}_{I^k_i}(t)dt
\end{equation}
for every Borel subset $B_k\subset[0,\oT_k]$. It follows from the imposed assumptions and Theorem~\ref{Thr3} that $\ox^k(t)$ and $\ox(t)$ are uniformly Lipschitzian on $[0,\oT]$, and thus by \eqref{e5.5} we have a constant $L$ such that
\begin{equation}\label{bdth}
\|\xi_{iy}^k\|\leq 2h^kL.	
\end{equation}
Set $\Lambda_i^k:=-\frac{\mu^k_0}{h^k}\xi_{iy}^k+p_{i+1}^k$ and $\Lambda^k(t):=\sum^{k-1}_{i=0}\Lambda_i^k\mathbbm{1}_{I^i_k}(t)$ with $I^i_k=[t_i^k,t_{i+1}^k)$, $i=0,\ldots,k-1$, and observe that 
\begin{equation}\label{bdL} 
\|\Lambda^k\|_{L^1}\le 2\mu^k_0 TL+h^k\disp\sum_{i=1}^k\|p^k_i\|.
\end{equation}
By the construction of $\Lambda^k(t)$, the discrete adjoint system \eqref{conx'} can be written as 
\begin{equation}\label{disadj}
\begin{aligned}
\dot{q}^k(t)&=-\nabla_x g\big(\ox^k(\nu^k(t)),\ou^k(\nu^k(t))\big)^*\Lambda^k(t)+\disp\sum_{j\in I_0(t,\Lambda^k(t)) }
\gg^k_j(t)x^j_\ast(t)+\disp\sum_{j\in I_>(t,\Lambda^k(t))}
\gg^k_j(t)x^j_\ast(t).
\end{aligned}
\end{equation}
It follows from (H3) and the discussions above that there is  $K>0$ such that $\|\nabla_x g(x,u)\|\leq K$ for any $x$ in a neighborhood of $\ox(t)$ as $t\in [0,\oT]\}$ and any $u\in U$. Since all the expressions in the statement of Theorem~\ref{Thr5} are positively homogeneous of degree one with respect to
$(\mu^k_0,p^k,q^k,\gg^k,\psi^k)$, the enhanced nontriviality condition \eqref{entc} allows us to normalize them by imposing the sequential equality
\begin{equation}\label{non-tri}
\begin{array}{ll}
\disp\mu^k_0+\|q^k(0)\|+\|p^k(\oT_k)\|&+\disp\sum_0^{k-1}\|\psi_i^k\|+\disp h^k\disp\sum_{\ell=0}^{k-1}\bigg\|\sum_{j\in I^{\ell+1}_{0k}(\Lambda^k_{\ell+1})}\gg_{\ell+1,j}^k x^j_{\ast}(t^k_{\ell+1})\bigg\|\\
&\disp+h^k\sum_{i=0}^{k-1}\sum_{j\in I^i_{>k}(\Lambda^k_i)}\gg^k_{ij}=1,
\end{array}
\end{equation}
where $\gg^k_{ij}\geq 0$ for all $j\in I^i_{>k}(\Lambda^k_i)$ according to \eqref{93}.

Without loss of generality, suppose that $\mu^k_0\to\mu$ as $k\to\infty$ for some $\mu\ge 0$ due to \eqref{non-tri}. To show the uniform boundedness of
$\{p^k_0,\ldots,p^k_k\}_{k\in\N}$ for all $i=0,\ldots, k-1$, $k\in\mathbb{N}$, we start with the observation from \eqref{conx} that
\begin{equation*}
p^k_{i+1}=p^k_i-h^k\nabla_x g(\ox^k_i,\ou^k_i)^*\Big(-\frac{1}{h^k}\mu^k_0\xi^k_{iy}+p^k_{i+1}\Big)+
h^k\sum_{j=1}^{s}\gg^k_{ij} x^j_\ast(t^k_i)
\end{equation*}
for all $i=0,\ldots,k-1$. By the above uniform boundedness $\|\nabla_x g(x,u)\|\leq K$, we get the estimate
\begin{equation}
\begin{aligned}\label{bdp0p1}
\|p^k_i\|\leq \(1+h^k K\)^{k-i}\|p^{k}_k\|+\sum_{\ell=i}^{k-1}\(1+h^k K\)^{\ell-i}K\mu^k_0\|\xi^k_{\ell y}\|+\sum_{\ell=i}^{k-1}\(1+h^k  K\)^{\ell-i}h^k\Big\|\sum_{j=1}^s\gg^k_{\ell j}x^j_\ast(t^k_i)\Big\|.
\end{aligned}
\end{equation}
On the other hand, it follows from \eqref{non-tri} that
\begin{equation}\label{2ndterm}
\sum_{\ell=i}^{k-1}\(1+h^k  K\)^{\ell-i}h^k\Big\|\sum_{j=1}^s\gg^k_{\ell j}x^j_\ast(t^k_i)\Big\|
\le e^{KT_k}\;\mbox{ for all }\;i=0,\ldots,k-1.
\end{equation}
Considering further the numbers
\begin{equation*}
A^k_i:=\sum_{\ell=i}^{k-1}\(1+h^k K\)^{\ell-i}K\mu^k_0\|\xi^k_{\ell y}\|
\end{equation*}
for $i=0,\ldots,k-1$, we deduce from \eqref{bdth} that
\begin{equation}
\begin{aligned}
A^k_i&\leq \sum_{\ell=i}^{k-1}\big(1+h^k K\big)^{\ell - i}K\mu^k_02h^kL
\leq  \mu^k_0 2 K L e^{KT_k}kh^k
\leq 2\mu^k_0KLT_k e^{KT_k},
\end{aligned}
\end{equation}
where the last term is bounded due to \eqref{non-tri} and the convergence of $T_k $ to $\oT$ as $k \to \infty$. Then recalling that the vectors $p^k_k=p_k(T_k)$ are uniformly bounded by \eqref{non-tri}, we obtain for $i=1, \ldots, k-1$ that
\begin{eqnarray*}
\|p^k_{i}\|&\le\big(1+Kh^k\big)^{k-i}\|p^{k}_k\|+A^k_i+ \sum_{\ell=i}^{k-1}(1+Kh^k)^{\ell-i}h^k\|\sum^s_{j-1}\gg^k_{\ell i}x^j_{\ast}(t^k_i)\|\le 2e^{KT_k }(1+\mu^k_0KLT_k)\le C
\end{eqnarray*}
with a suitable constant $C$. Thus the boundedness of the sequence $\{p^k_{0}\}$ follows from \eqref{bdp0p1} and the boundedness of $\{p^k_{i}\}_{1\le i\le k}$,
which therefore justifies the boundedness of the entire bundle $\{( p^k_{0},\ldots,p^k_k)\}_{k\in\N}$. Taking into account that the subgradient sets for locally Lipschitzian functions are bounded and employing \eqref{c:6.14'} together with
\eqref{non-tri} tell us that the sets $\partial\la -\frac{\mu^k_0\xi^k_{iy}}{h^k}+p^k_{i+1},-g\ra (\ox^k_i,\ou^k_i)$ in \eqref{conx} are uniformly bounded for all $i=0,\ldots,k-1$. Hence there exists a constant $M>0$ ensuring that
$$
\Bigg\|\Bigg(\frac{p^k_{i+1}-p^k_i}{h^k} -\sum_{j\in I^i_{0k}(-\frac{\mu^k_0\xi^k_{iy}}{h^k}+p^k_{i+1})
\cup I^i_{>k}(-\frac{\mu^k_0\xi^k_{iy}}{h^k}+p^k_{i+1})}\gg^k_{ij} x^j_*(t^k_i),-\frac{\mu^k_0\xi^k_{iu}}{h^k}-\frac{\psi^k_i}{h^k}\Bigg)\Bigg\|\leq M.
$$
This allows us to deduce from \eqref{c:6.25} the estimate 
\begin{equation}\label{psik} 
\|\psi^k(t)\|\leq M+\mu^k_0\|\xi^k_u(t)\|\;\mbox{ for }\;t\in [t^k_i,t^k_{i+1}).
\end{equation}
By \eqref{c:6.14'} and \eqref{psik}, we get the boundedness of $\{\psi^k(\cdot)\}$ in $L^2([0,\oT_k];\R^d)$, and so there exists a subsequence of $\{\psi^k(\cdot)\}$, which weakly converges to some function $\psi(\cdot)\in L^2([0,\oT];\R^d)$ as $k\to \infty$.\vspace*{-0.05in} 

Next we show that the functions $q^k(\cdot)$ have uniformly bounded variations, i.e., the norm sequence
$\{\|\dot{q}^k\|_{L^1}\}$ is uniformly bounded. It follows from the definitions of $I^i_{0k}$ and $I^i_{>k}$ in \eqref{c56} that if $\ell\in I^i_{0k}(\Lambda^k_i)$, then 
$$
\la x^{\ell}_{\ast}(t^k_i),p^k_{i+1}\ra=\frac{\mu^k_0}{h^k}\la x^{\ell}_{\ast}(t^k_i),\xi_{iy}^k\ra,
$$
due to the imposed LICQ, and so we have the uniform estimate 
\begin{equation*}
\bigg\la x^{\ell}_{\ast}(t^k_i), \dfrac{p^k_{i+1}-p_i^k}{h^k}\bigg\ra
\leq \dfrac{\mu^k_0}{(h^k)^2}\la x^{\ell}_{\ast}(t^k_i), \xi_{(i)y}^k-\xi_{(i-1)y}^k\ra,
\end{equation*}
where the equality holds if $\ell\in I^{i-1}_{0k}(\Lambda^k_{i-1})$ and the strict inequality holds if $\ell\in I^{i-1}_{>k}(\Lambda^k_{i-1})$. Similarly, 
\begin{equation*}
\bigg\la x^{\ell}_{\ast}(t^k_i), \dfrac{p^k_{i}-p_{i+1}^k}{h^k}\bigg\ra
\geq \dfrac{\mu^k_0}{(h^k)^2}\la x^{\ell}_{\ast}(t^k_i), \xi_{(i-1)y}^k-\xi_{iy}^k\ra.
\end{equation*}
Therefore, for all $ \ell\in I^i_{0k}(\Lambda^k_{i})$ as $i \in \{1, \ldots, k-1\}$ and $k\in \N$, we obtain that 
\begin{equation}
\label{adjI0}
\Bigg|\bigg\la x^{\ell}_{\ast}(t^k_i), \dfrac{p^k_{i+1}-p_i^k}{h^k}\bigg\ra\Bigg|
\leq \dfrac{\mu^k_0}{(h^k)^2}|\la x^{\ell}_{\ast}(t^k_i), \xi_{(i)y}^k-\xi_{(i-1)y}^k\ra|.
\end{equation}
Combining \eqref{adjI0} and \eqref{disadj} and then summing  them over $\ell\in I^i_{0k}(\Lambda^k_{i})$ yield
\begin{equation*}
\begin{aligned}
\sum_{\ell\in I^i_{0k}(\Lambda^k_{i})} \Bigg|\sum_{j\in I^i_{0k}(\Lambda^k_{i})}\gg^k_{ij}\la x^j_\ast (t^k_i), x^{\ell}_{\ast}(t^k_i)\ra \Bigg|
&\le \dfrac{\mu^k_0}{(h^k)^2}\sum_{\ell\in I^i_{0k}(\Lambda^k_{i})}|\la x^{\ell}_{\ast}(t^k_i), \xi^k_{iy}-\xi^k_{(i-1)y}\ra|\\
&+\sum_{\ell\in I^i_{0k}(\Lambda^k_{i})}|\la x^{\ell}_{\ast}(t^k_i), \nabla_x g(\ox^k_i, \ou^k_i)^*\Lambda^k_ih^k\ra|+s\sum_{j\in I^i_{>k}(\Lambda^k_{i})}\gg^k_{ij}
\end{aligned}
\end{equation*}
for all $i \in \{1, \ldots, k-1\}$ and $k\in \N$. Observe further that the functions $\disp\sum _{i=1}^k\dfrac{\xi^k_{iy}-\xi^k_{(i-1)y}}{(h^k)^2}\mathbbm{1}_{I^i_k}(t)$ are bounded in $L^1$ in $k$. Indeed, it follows from Theorem~\ref{Thr1} that
$$
\Bigg\|\dfrac{\xi^k_{(i+1)y}-\xi^k_{iy}}{(h^k)^2}\Bigg\|=\frac{1}{h^k}\Bigg\|\dfrac{\bar{x}^k_{i+2}-\bar{x}^k_{i+1}}{h^k}
-\dfrac{\bar{x}(t^k_{i+2})-\bar{x}(t^k_{i+1})}{h^k} \Bigg\| +\frac{1}{h^k}\Bigg\|\dfrac{\bar{x}^k_{i+1}-\bar{x}^k_i}{h^k}
-\dfrac{\bar{x}(t^k_{i+1})-\bar{x}(t^k_i)}{h^k}\Bigg\|\leq \dfrac{4L}{h^k},
$$
and therefore we have the estimate
$$
\int_0^{\oT}\sum_{i=0}^{k-1}\bigg\|\dfrac{\xi^k_{(i+1)y}-\xi^k_{iy}}{(h^k)^2}\bigg\|\mathbbm{1}_{I^i_k}(t)dt\le 4L.
$$
The functions $\sum^{k-1}_{i=0}\nabla_x g(\ox^k_i, \ox^k_i)^*\Lambda^k_i\mathbbm{1}_{I^i_k}(\cdot)$ are bounded in $L^1$ uniformly in $k$ by \eqref{bdL}, while the functions
$\disp \sum^{k-1}_{i=0}\sum_{j\in I^i_{>k}(\Lambda^k_{i})}\gg^k_{ij}\mathbbm{1}_{(I^i_k)}(\cdot)$ are bounded in $L^1$ uniformly in $k$ due to the normalization condition \eqref{non-tri}. Similarly,
\begin{equation*}
\begin{aligned}
\Big|\la v, \sum_{j\in I^i_{0k}(\Lambda^k_{i})}\gg^k_{ij}x^j_\ast(t^k_i)\ra\Big|
&\le \dfrac{\mu^k_0}{(h^k)^2}\| v\|\cdot\| \xi^k_{iy}-\xi^k_{(i-1)y}\|+\|v\|\cdot\| \nabla_x g(\ox^k_i, \ou^k_i)^*\|\cdot\|\Lambda^k_i\| +s\sum_{j\in I^i_{>k}(\Lambda^k_{i})}\gg^k_{ij}
\end{aligned}
\end{equation*}
for each $v\in \span \{x^j_{\ast}(t^k_i)\;|\;j\in I^i_{0k}(\Lambda^k_i) \}=: V^k_i,$ $i=0,\ldots, k-1$.
Let now $v$ be any vector in $\R^n$. For each $k\in \N$ and $i=1,\ldots, k$, 
denote $v^k_i:=\proj_{V^k_i}(v)$. 
Then we have
\begin{equation}
\label{5.35}
\begin{aligned}
&h^k\disp\sum^{k}_{i=1}\Big|\la  v, \sum_{j\in I^i_{0k}(\Lambda^k_i)}\gg^k_{ij}x^j_{\ast}(t^k_i)\ra \Big|=\disp\sum^{k}_{i=1}\Big|\la v^k_i, \sum_{j\in I^i_{0k}(\Lambda^k_i)}\gg^k_{ij}x^j_{\ast}(t^k_i)\ra\Big|h^k \\
&\leq \disp\sum^{k}_{i=1}\dfrac{\mu^k_0}{h^k}\| v^k_i\|\cdot\| \xi^k_{iy}-\xi^k_{(i-1)y}\|+ h^k \disp\sum^{k}_{i=1}\| v^k_i\|\cdot\| \nabla_x g(\ox^k_i, \ou^k_i)\|\cdot\|\Lambda^k_i\|+2sh^k\sum^k_{i=1}\sum_{j\in I^i_{>k}(\Lambda^k_i)}\gg^k_{ij},
\end{aligned}	
\end{equation} 
where the right-hand side of the inequality is bounded uniformly in $k$ due to the above arguments and the normalization condition \eqref{non-tri}. For all $v\in \R^n$, we deduce from \eqref{disadj} that 
$$\la v, \dot q^k(t)\ra =\la v, \dfrac{p^k_{i+1}-p^k_i}{h^k}\ra =-\la v,\nabla_x g(\ox^k_i,\ou^k_i)\Lambda^k_i\ra + \la v^k_i, \sum_{j\in I^i_{0k}(\Lambda^k_i)}\gg^k_{ij}x^j_{\ast}(t^k_i)\ra +\la v, \sum_{j\in I^i_{>k}(\Lambda^k_i)}\gg^k_{ij}x^j_{\ast}(t^k_i)\ra
$$
whenever $t\in (t^k_i, t^k_{i+1})$. It follows from \eqref{5.35} that
$$
\begin{aligned}
\int^{T_k}_0|\la v, \dot q^k(t)\ra |dt\leq h^k \disp\sum^{k}_{i=1} \| v\|\cdot\|\nabla_x g(\ox^k_i, \ou^k_i)\|\cdot\|\Lambda^k_i\|+ h^k \disp\sum^{k}_{i=1}|\la v^k_i, \sum_{j\in I^i_{0k}(\Lambda^k_i)}\gg^k_{ij}x^j_{\ast}(t^k_i)\ra|+h^k\disp\sum^{k}_{i=1}|\la v, \sum_{j\in I^i_{>k}(\Lambda^k_i)}\gg^k_{ij}x^j_{\ast}(t^k_i)\ra|
\end{aligned}
$$
for all vectors $v\in\R^n$, and that the right-hand side of the above inequality is bounded uniformly in $k$ and $v$. This ensures that the sequence $\{\|\dot q^k\|_{L^1([0,\oT_k];\R^n)}\;|\;k\in \N\}$ is bounded. Applying now Helly's selection theorem gives us a function of bounded variation $q(\cdot)$
such that $q^k(t)\to q(t)$ as $k \to\infty$ pointwise on $[0,\oT]$.\vspace*{-0.05in}

Observe that since the sequence of measures $\{dq^k:=\dot q^k dt\}$ has uniformly bounded total variation, its subsequence converges weakly$^\ast$ in $C^*([0,\oT];\R^s)$ to a measure $dq$. It follows from \eqref{c:59}, \eqref{non-tri}, and the uniform boundedness of $q^k(\cdot)$ on $[0,\oT_k]$ that the sequence $\{p^k(\cdot)\}$ is bounded in $W^{1,2}([0,\oT_k];\R^{n})$ and thus weakly compact in this reflexive space. Moreover, the normalization condition \eqref{non-tri} ensures the positive measures $d\gg^k_>:= \disp\sum^{k-1}_{i=0}\sum_{j\in I^i_{>k}(\Lambda^k_i)}\gg^k_{ij}\mathbbm{1}_{I^k_i}dt$ have uniformly bounded total variations. Therefore, it follows from \eqref{c:6.30} that the measures $d\gg^k_0:=\disp\sum^{k-1}_{i=0}\sum_{j\in I^i_{0k}(\Lambda^k_i)}\gg^k_{ij}\mathbbm{1}_{I^k_i}dt$ have uniformly bounded variations as well. Since $\gg^k_{ij}=0$ for all $j\notin I^i_{0k}(\Lambda^k_i)\cup I^i_{>k}(\Lambda^k_i)$, we get, up to a subsequence, that $\gg^k\st{w^*}{\to}\gg$, $\gg^k_>\st{w^*}{\to}\gg_>$, and $\gg^k_0\st{w^*}{\to}\gg_0$ with $\gg=\gg_>+\gg_0$, where $w^*$ stands for the weak$^*$ convergence in the corresponding space. Mazur's theorem tells us that there exists a sequence of convex combinations of $\{\dot p^k(\cdot),\psi^k(\cdot)\}$, which converges to some $(\dot p(\cdot),\psi(\cdot))\in L^2([0,\oT];\R^n)\times L^2([0,\oT];\R^n)$ a.e.\ pointwise on $[0,\oT]$. This gives us \eqref{c:6.6} by passing to the limit along \eqref{conx'} as $k\to\infty$ with the usage of \eqref{c:6.14} and \eqref{c:6.14'} up to choosing the right continuous representative of $q$. We also get the convergence
\begin{equation*}
\Big\|\int_t^{\oT_k}\sum_{j=1}^s\gg^k_j(\tau)x^j_\ast (\tau)d\tau-\int_{(t,\oT]}\sum_{j=1}^sd\gg^j(\tau)x^j_{\ast}(\tau)\Big\|\to 0\;
\textrm{ as }\;k\to\infty
\end{equation*}
for all $t\in[0,\oT_k]$ except a countable subset of $[0,\oT_k]$ by the weak$^*$ convergence of the measures $\gg^k$ to $\gg$ in
$C^*([0,\oT];\R^n)$; cf.\  \cite[p.\ 325]{v} for similar arguments. Hence we have the convergence
\begin{equation*}
\int^{\oT_k}_t\sum_{j=1}^s\gg^k_j(\tau)x^j_\ast (\tau)d\tau\to\int_{(t,\oT]}\sum_{j=1}^sd\gg^j(\tau)x^j_{\ast}(\tau)\;\textrm{ on }\;[0,\oT]\;
\textrm{ as }\;k\to\infty
\end{equation*}
and thus arrive at \eqref{c:6.9} by passing to the limit in \eqref{c:6.29}. The claimed condition $p(\oT)=q(\oT)$ follows directly by passing to the limit in the equalities $p^k(\oT_k)=q^k(\oT_k)$, $k\in\N$. The second complementary implication in \eqref{41} follows from \eqref{96} under LICQ while arguing by contradiction with the usage of the established a.e.\ pointwise convergence
of the functions involved therein. To verify now \eqref{co}, recall from Theorem~\ref{Thr5} that $\psi^k_i\in N(\ou^k_i;U)$ for $i=0,\ldots,k-1$. By the construction of $\psi^k(\cdot)$ in \eqref{c:6.25}, the piecewise constant extension of $\ou^k_i$ to $[0,\oT_k]$, and the conic structure of $N(\cdot;U)$ we get that
\begin{equation}\label{coN}
\psi^k(t)\in N\big(\ou^k(t);U\big)\;\textrm{ for all }\;t\in[t^k_i,t^k_{i+1})\;\textrm{ and }\;i=0,\ldots,k-1.
\end{equation}
Then the desired result follows by passing to the pointwise limit in \eqref{coN} along a subsequence of $k\to\infty$ with employing the strong $L^2$-convergence of $\ou^k(\cdot)\to\ou(\cdot)$ from 
Theorem~\ref{Thr1}, the robustness of the normal cone with respect to perturbations of the initial point, the strong $L^2$-convergence of convex combinations of $\psi^k(\cdot)$ to $\psi(\cdot)$, and the boundedness of $\psi(\cdot)$ on $[0,\oT]$ due to \eqref{c:6.6} under (H3) and $q(\cdot)\in BV([0,\oT];\R^n)$.\vspace*{-0.05in}

Finally in this step, it remains to verify the fulfillment of the tangential maximization condition in \eqref{c:6.6'} and its global version in \eqref{max}. Note that the duality correspondence in \eqref{dua} generated by any tangent set $T(\ou(t);U)$ always yields the convexity of $N(\ou(t);\O)$. We deduce \eqref{c:6.6'} directly from \eqref{co}, \eqref{dua}, and the
inclusions $\psi^k_i\in N(\ou^k_i;U)$ for $i=0,\ldots,k-1$ of Theorem~\ref{Thr6} as $k\to\infty$. If $U$ is convex, the maximization condition \eqref{max} follows from \eqref{c:6.6'}
due to the structure \eqref{nc} of the normal cone in convex analysis.\\[1ex]
{\bf Step~3:} {\em Transversality and complementary slackness at the optimal final time.} Let us verify verify the endpoint inclusion \eqref{42}, which combines the transversality
conditions on $p(\oT)$ with the additional condition on the optimal
time interval $[0,\oT]$. First we examine the passage to limit on the left-hand side of \eqref{p_kk} as $k\rightarrow\infty$. Having
\begin{eqnarray*}
(-p^{k}(\oT_{k}),\;2\mu^k_{0}(\oT-\oT_{k}))\rightarrow(-p(\oT),0)
\end{eqnarray*}
and considering $\varrho_{k}$ from \eqref{vrho} give us the relationships
\begin{eqnarray*}
\begin{array}{ll}
\varrho_k:=\disp\sum_{i=0}^{k-1}\left[\dfrac{i}{k}\n\(\frac{\ox^k_{i+1}-\ox^k_i}{h^k}-
\displaystyle\dot{\ox}(t_i),\ou^k_i-\ou(t_i)\)\en^2-\dfrac{i+1}{k}\n\(\frac{\ox^k_{i+1}-\ox^k_i}{h^k}-
\dot{\ox}(t_{i+1}),\ou^k_i-\ou(t_{i+1})\)\en^2\right]\\
=\disp\mathbbm{1}_{I^i_k}(\cdot)\mathbbm{1}_{I^i_k}(\cdot)\frac{1}{\oT_{k}}\sum_{i=0}^{k-1}\[i\int_{t_{i}}^{t_{i+1}}\n\(\dot{i
\bar{x}}^{k}(t)-\dot{\bar{x}}(t_{i}),\ou^k_i-\ou(t_i)\)\en^{2}dt-(i+1)\int_{t_{i}}^{t_{i+1}}\n\(\dot{
\bar{x}}^{k}(t)-\dot{\bar{x}}(t_{i+1}),\ou^k_i-\ou(t_{i+1})\)\en^{2}dt\]\\
\disp\sim\frac{1}{\oT_{k}}\sum_{i=0}^{k-1}\int_{t_{i}}^{t_{i+1}}\n\(\dot{\bar{x}}^{k}(
t)-\dot{\bar{x}}(t),\ou^k(t)-\ou(t)\)\en^{2}dt=\frac{1}{\oT_{k}} \int_{0}^{\oT_{k}}\n\(\dot{\bar
{x}}^{k}(t)-\dot{\bar{x}}(t),\ou^k(t)-\ou(t)\)\en^{2}dt\rightarrow 0
\end{array}
\end{eqnarray*}
as $k\to\infty$ due to the strong convergence in Theorem~\ref{Thr3}.
To evaluate $\bar{H}^{k}$ in \eqref{Hk}, we employ \eqref{vel}, \eqref{t}, and the
convergences in Theorem~\ref{Thr3} together with the uniform convergence $p^{k}(\cdot)\rightarrow p(\cdot)$ on $[0,\oT]$ while getting
\begin{eqnarray*}
\bar{H}^k&:=&\dfrac{1}{k}\sum_{i=0}^{k-1}\la p^k_{i+1},y^k_i\ra =\frac{1}{\oT_{k}}\sum_{i=0}^{k-1}\int_{t_{i}}^{t_{i+1}}\langle
p^{k}(t_{i+1}),\dot{\bar{x}}^{k}(t)\rangle dt\\
&\sim& \frac{1}{\oT_{k}}\sum_{i=0}^{k-1}\int_{t_{i}}^{t_{i+1}}\langle p^{k}(t),
\dot{\bar{x}}^{k}(t)\rangle dt=\frac{1}{\oT_{k}}\int_{0}^{\oT_{k}}\langle p^{k}(t),\dot{\bar{x}}^{k}
(t)\rangle dt\\
&\rightarrow&\frac{1}{\bar
{T}}\int_{0}^{\oT}\langle p(t),\dot{\bar{x}}(t)\rangle dt:=\bar{H}\;\;\mbox{as}\;\;k\rightarrow\infty
\end{eqnarray*}
by the Lebesgue dominated convergence theorem.\vspace*{-0.05in}

To proceed further, we examine the passage to the limit on the right-hand side expression in \eqref{p_kk}. Note that 
\begin{eqnarray*}
\Limsup_{k\rightarrow\infty}\partial(\mu^k_0\varphi)(\bar{x}^k(\oT_k),\oT_k)=\mu\partial\varphi(\bar{x}(\oT),\oT)
\end{eqnarray*}
due to the robustness of the basic subdifferential \eqref{sub1}. Relying on the discrete necessary optimality conditions of Theorem~\ref{Thr6}, define $\eta^k(\oT):=\eta^{k}_k$ and deduce from the normalization of the nontriviality conditions in \eqref{ntc} that a subsequence of $\{\eta^{k}_k\}$ converges, without relabeling, to some vector $(\eta^{1}(\oT),\ldots,\eta^{k}(\oT))$. It follows from \eqref{p_kk} and representation \eqref{F} that for each $k\in\N$ we have the inclusion
\begin{equation}\label{c:72}
-p^k_k-\mu^k_0\vt^k_k= \sum_{j\in I^k_k(\ox^{k}_k)}\eta^k_{kj} x^{j}_\ast(t^k_k)
\in N(\ox^{k}_k;C^k_k),
\end{equation}
where $\eta^k_{kj}=0$ if $j\in\{1,\ldots,s\}\setminus I^k_k(\ox^{k}_k)$. Denoting $\zeta^k:=
\sum_{j\in I^k_k(\ox^k_k)}\eta^k_{kj} x^k_\ast(t^k_k)$ and using the boundedness of $\mu^k_0$ by \eqref{non-tri} together with the convergence of $\{p^k_k\}$ and $\{\ox^k_k\}$ allows us to select a subsequence of $\{\zeta^k\}$ converging to  some $\zeta\in\R^n$. It follows from the robustness of the normal cone in \eqref{c:72}, the convergence of $\ox^k_k\to \ox(\oT)$, and
the inclusion $I^k_k(\ox^k_k)\subset I(\oT,\ox(\oT))$ for all $k$ sufficiently large that $\zeta\in N(\ox(\oT);C(\oT))$. Observe that the set $\Omega^k_x\times \Omega_T^k$ in \eqref{Omega} admits the representation
$$
\Omega^k_x\times \Omega_T^k=\big\{ (x^k_k,T_k)\in R^{n}\times[0,\infty)\;\big|\;\text{dist}((x^k_k,T_k) ,(\O_x\times \O_T))\leq \delta^k\big\}.
$$
Applying \cite[Proposition~2.7]{m88} if $(\ox^k_k,\oT_k)\in \O_x\times \O_T$ and the normal cone definition if $(\ox^k_k,\oT_k)\notin \O_x\times \O_T$, we get
$$\Limsup_{k\rightarrow\infty}N((\ox^k_k,\oT_k);\Omega^k_x\times \Omega_T^k)=N((\ox(\oT),\oT); \O_x\times \O_T).$$
Passing now to the limit in \eqref{p_kk} as $k\to\infty$ verifies the transversality inclusion 
$$
(-p(\oT)-\sum_{j\in I(\oT,\ox(\oT))}\eta_j(\oT) x^{j}_\ast(\oT),\bar{H})\in \mu\partial\varphi(\bar{x}(\oT),\oT)+N((\ox(\oT),\oT); \O_x\times \O_T),
$$
together with the implication in \eqref{42}. Finally, the fulfillment of the endpoint complementary slackness conditions in \eqref{41a} follows directly from the above proof by passing to the limit as $k\to\infty$ in their discrete counterparts established in \eqref{eta1} of Theorem~\ref{Thr6}.\\[1ex]
{\bf Step~4:} {\em Measure nonatomicity.} Take $t\in[0,\oT)$ with $\la x^j_\ast(t),\ox(t)\ra<c_j(t)$ for all $j=1,\ldots,s$ and by continuity of $\ox (\cdot)$ find a neighborhood $V_t$ of $t$ such that $\la x^j_\ast(\tau),\ox(\tau)\ra<c_j(\tau)$ whenever
$\tau\in V_t$ and $j=1,\ldots,s$.
Invoking Theorem~\ref{Thr3} tells us that $\la x^j_\ast(t^k_i),\ox^k(t^k_i)\ra<c_j(t^k_i)$ if $t^k_i\in V_t$ for all $j=1,\ldots,s$ and
$k\in\N$ sufficiently large. Then we deduce from \eqref{94} that $\gg^k_j(t)=0$ as $j=1,\ldots,s$ on any Borel subset $V$ of $V_t$. Hence
\begin{equation}\label{nonatom}
\|\gg^k_j\|(V)=\disp\int_Vd\|\gg^k_j\|=\int_V\|\gg^k_j(t)\|dt=0
\end{equation}
by the construction of $\gg^k$ in \eqref{c:6.34}. Passing to the limit therein and taking into account the measure convergence obtained
in Step~2, we get $\|\gg_j\|(V_t)=0$, which justifies the claimed nonatomicity condition.\\[1ex]
{\bf Step~5:} {\em Nontriviality and support conditions.} We begin with the proof of the nontriviality condition \eqref{e:83} under the general assumptions of the theorem. Defining $\eta^k_j(t):=\sum_{i=0}^{k-1}\eta^k_{ij}\mathbbm{1}_{[t^k_i,t^k_{i+1})}(t)$ and remembering that $\bar{x}^k\to\bar{x}$
in $W^{1,2}$ as $k\to\infty$, we obtain that $\eta^k_j\to\eta^j$ strongly in $L^2$. If the normal cone is active along $\bar{x}(t)$ only for some $t\in E\subset [0,\oT]$ with nonempty interior, then Definition~\ref{active cone} of the active normal cone tells us that for all $ t\in \mbox{int}(E_0)$ there exists $k(t)$ such that $\eta^k_j(t)>0$ whenever $k\ge k(t)$ and $j\in I(t,\ox)$. Consequently, for each $t\in \mbox{int}(E_0)$ and $j\in I(t,\ox)$, we get by recalling condition \eqref{96} obtained under LICQ that
\begin{equation*}
\Big\la x^j_\ast(t^k_i),-\frac{\mu^k_0\xi^k_{iy}}{h^k}+p^k_{i+1}\Big\ra=0.
\end{equation*}
The latter equality being combined with condition \eqref{c56} implies that $j\in I^i_{0k}(\Lambda^k_{i})$, and so that $I_{>k}(\Lambda^k_i)=\emptyset$ for all $k\geq k(t)$. 
This allows us to check that the support of the weak$^{\ast}$-limit of the last summand in \eqref{disadj} doesn't intersect the interior of $E_0$. Thus we arrive at the claimed support condition \eqref{suppcond}.\vspace*{-0.05in}

To proceed with the proof of the general nontriviality condition \eqref{e:83}, suppose on the contrary that $\mu=0$ and $p(t)=0$ for all $t\in[0,T]$ together with $\gg_0=\gg_>=0$. This implies by recalling \eqref{c:6.9} that $q=0$. Then $\mu^k_0\to 0$, $p^k\to 0$ uniformly, and $q^k\to 0$ pointwise. Moreover, we actually have from the above that $q^k\to 0$ in $L^1$ as $k\to\infty$. Since $\gg^k_>$ are positive measures, it follows that $\disp \lim_{k\to \infty}h^k\sum^{k-1}_{i=0}\sum_{j\in I^i_{>k}(\Lambda^k_i)}\gg^k_{ij}=0$. As a consequence, we get by recalling \eqref{bdL} that all the summands on the right-hand side of \eqref{5.35} vanish as $k\to \infty$ uniformly with respect to the unit vectors $v$. This tells us in turn that 
$$
\disp\lim_{k\to \infty}h^k\sum^{k-1}_{\ell=0}\Big\|\sum_{j\in I^{\ell+1}_{0k}(\Lambda^k_{\ell+1})}\gg^k_{\ell+1, j}x^j_{\ast}(t^k_{\ell+1})\Big\|=0.
$$ 
The above arguments together with \eqref{conx} ensure also that $\lim_{k\to \infty}\sum^{k-1}_{i=0}\|\psi^k_i\|=0$, which clearly contradicts \eqref{non-tri} and therefore justifies \eqref{e:83}.\vspace*{-0.05in}

It remains to verify the fulfillment of the {\em enhanced nontriviality} in the theorem under the imposed interiority assumption. Suppose on the contrary that $\mu=0$ for all $t\in[0,\oT]$ while $\la x^j_{\ast}(t),\ox(t)\ra<c_j(t)$ for such $t\in[0,\oT]$
and all $j=1,\ldots,s$. It follows from the nonatomicity condition together with 
\eqref{c:6.9} that $\gg_>=\gg_0=0$, and so
\begin{equation}\label{qA}
q(t)=p(t)-\int_{(t,\oT]}\sum_{j=1}^sd\gg^j(\tau)x^j_{\ast}(\tau)d\tau=p(t)\;\textrm{ for all }\;t\in[0,\oT]\setminus A,
\end{equation}
where $A\subset[0,\oT]$ is a countable set. Moreover, it follows from \eqref{41} that $\eta(t)=0$ on $[0,\oT]$, which being combined with the first inclusion in \eqref{42} gives us $p(\oT)=0$. Since $p$ is a solution to a linear homogeneous equation with the vanishing final condition, it must be identically $0$, and thus the failure of the enhanced transversality contradicts the fulfillment of the general nontriviality condition \eqref{e:83}. This completes the proof of the theorem. $\h$\vspace*{-0.1in}

\begin{remark}\label{rem-compar} Let us briefly compare the necessary optimality conditions for controlled sweeping processes established in Theorem~\ref{Thr7} with the recent results in this direction related to our paper:\\[1ex]
$\bullet$ The results of Theorem~\ref{Thr7} significantly improve the necessary optimality conditions obtained in \cite{cmn18a} for a particular case of problem $(P)$, where the duration of the process is fixed, the endpoint constraints are absent, and the sweeping set $C$ is not moving. Even in this particular case, our novel results provides essentially new information on local optimal solutions to $(P)$ including the new support condition in both active and nonactive settings for the normal cone along the optimal solution (this important notion was missed in \cite{cmn18a}). Furthermore, the nontriviality condition in \cite[Theorem~7.1]{cmn18a}, which is subject to a possibly worse degeneration, requires LICQ. Now we establish the new nondegeneracy condition involving both positive and signed measures, the new enhanced nontriviality condition obtained in the normal (Fritz John) form, consider the general case of measurable optimal controls, and avoid the smoothness assumption on the mapping $g(x,u)$ in \eqref{Problem} with respect to both variables $(x,u)$ as well as the additional full rank assumption on $\nabla_u g$ along the solution. This is done by developing the novel version of the method of discrete approximations of its independent interest.\\[1ex]
$\bullet$ Observe that if the normal cone is inactive along the optimal couple $(\ox,\ou)$, i.e., in \eqref{37} all coefficients $\eta^j$, $j=1,\ldots,s$, vanish on $[0,\oT]$, then $(\ox,\ou)$ is an optimal couple also for the (classical) dynamic $\dot{\ox}=g(x,u)$, subject to the state constraint $x(t)\in C(t)$ for all $t\in[0,\oT]$. Nevertheless, the index set $I_0$ may be nonempty, and so the signed measured $\gg_0$ may be nonvanishing. Therefore, in this case, the necessary conditions in classical state constrained optimal control problems are not recovered. This is due to the fact that the class of admissible trajectories is larger than in the classical state constrained problems as it may contain competitors of $(\ox,\ou)$ for which the normal cone is active on a set of positive measure.\\[1ex]
$\bullet$ In \cite{her-pall}, Hermosilla and Palladino address a fixed-time version of problem $(P)$ without additional endpoint constraints under a set of assumptions that are generally different from ours. They developed a powerful continuous-time penalization technique and derive necessary optimality conditions that are essentially different from those in Theorem~\ref{Thr7}. Note that the necessary optimality conditions in \cite{her-pall} contain a full sweeping counterpart of the Pontryagin maximum principle, while with a nontriviality condition that may degenerate and thus requires further investigations. The possible degeneration, like in our case, is due to the presence of signed measures.\\[1ex]
$\bullet$  De Pinho, Ferreira and Smirnov have recently considered in \cite{pfs23} a fixed-time sweeping optimal control problem similar to that in \cite{her-pall}, but with endpoint constraints and under a different set of assumptions in comparison with \cite{her-pall} and our current paper. In particular, the assumptions of \cite{pfs23} cover the moving polyhedral sets $C(t)$ in the sweeping process \eqref{Problem} only under certain relationships between the generating vectors of the polyhedron, which may fail to hold even for polyhedral sets in $\R^2$. In \cite{pfs}, the same authors develop a continuous-time smooth approximation method allowing them to derive a sweeping counterpart of the Pontryagin maximum principle by passing to the limit from optimal control of smooth ODE systems that are constraint free and thus leads to a normal form of the maximum principle. The necessary optimality conditions obtained in \cite{pfs,pfs23} are independent from those in \cite{her-pall} and in our paper. See also the new paper by  Nour and Zeidan \cite{vera} for further developments.\\[1ex]
$\bullet$ In \cite{zeidan}, the authors develop smooth approximation procedures and derive necessary optimality conditions of the maximal principle type for sweeping-related optimal control problems governed by subdifferential operators on fixed-time intervals, where the imposed assumptions are generally different from those discussed above.
\end{remark}\vspace*{-0.2in}

\section{Numerical Example}\label{examples}\vspace*{-0.05in}

In this section, we elaborate in detail a two-dimensional example of a free-time sweeping optimal control problem with endpoint constraints. The given example illustrates how the obtained necessary optimality condition work and allow us to explicitly determine optimal solutions.\vspace*{-0.1in}
 
\begin{example}\label{Ex-1}
In $\R^2$, consider problem (SP) formulated as follows:
\begin{equation*}
\mbox{minimize }\;J[x,u,T]=\varphi(x_1(T),T):=T+\disp\frac{1}{2}(x_1(T)-\alpha)^2
\end{equation*}
subject to the sweeping dynamics and endpoint constraints
\begin{equation}\label{dyn_ex8.2}
\left\{\begin{matrix}
\begin{pmatrix}
\dot{x}_1\\
\dot{x}_2
\end{pmatrix}=-N\left(\begin{pmatrix}
x_1\\
x_2
\end{pmatrix};C(t)\right)+\begin{pmatrix}
0\\
u
\end{pmatrix}\\
\textrm{with }\;\begin{pmatrix}
x_1\\
x_2
\end{pmatrix}(0)=\begin{pmatrix}
0\\
0
\end{pmatrix},\; x_2(T)= 1,
\end{matrix}\right.
\end{equation}
where $\alpha\in \R$ is a parameter, $g(x,u):=(0,u), \,u\in U:= [-2,2]$ for $t\in [0,T]$, and the moving set is defined by
\begin{equation}\label{bis}
C(t):=\Big\{x=(x_1,x_2)\;\Big|\;\Big\la(x_1,x_2),\Big(\frac{1}{\sqrt{2}},\frac{1}{\sqrt{2}}\Big)\Big\ra\le\dfrac{1}{\sqrt{2}}-\dfrac{t}{\sqrt{2}}\Big\}.
\end{equation}
It follows from \cite{CCMN21}, which proof can be easily modified for free-time sweeping control systems, that the formulated problem $(SP)$ admits an optimal solution $(\ox(\cdot),\ou(\cdot),\oT)$, where $(\ox,\ou)\in W^{1,2}\times L^2$. To proceed with the usage of the obtained necessary optimality conditions to explicitly calculate optimal solutions to (SP), observe first that we can represent the dynamics in the form
\begin{eqnarray}\label{dyna}
\begin{pmatrix}
\dot{x}_1\\
\dot{x}_2
\end{pmatrix}(t)=
\begin{pmatrix}
0\\
u(t)
\end{pmatrix}
-\eta(t)\begin{pmatrix}
\frac{1}{\sqrt{2}}\\[1ex]
\frac{1}{\sqrt{2}}
\end{pmatrix},\;\mbox{ where }\;\eta(t)\ge 0\;\textrm{ for a.e. }\;t\in[0,\oT]
\end{eqnarray}
accompanied by the relationships
$$
\Big\la (x_1(t), x_2(t)),\Big(\dfrac{1}{\sqrt{2}},\dfrac{1}{\sqrt{2}}\Big)\Big\ra<\dfrac{1}{\sqrt{2}}-\dfrac{t}{\sqrt{2}}\Longrightarrow \eta(t)=0.$$
Employing Theorem~\ref{Thr7} along a $W^{1,2}\times L^2\times\R_+$-local minimizer $(\ox(\cdot),\ou(\cdot),\oT)$, we get:

{\bf 1.} $(\dot p(t), \psi(t))\in{\rm co}\,\partial \( q_2(t)\ou(t)\)(\ox(t), \ou(t)),$ where $q(t)=(q_1(t), q_2(t))$ for a.e. $t\in[0,\oT]$ (by \eqref{c:6.6}).

{\bf 2.} $\psi(t)\in N\(\ou(t) ; U\) \;\textrm{for a.e. }\;t\in[0,\oT]$ (by \eqref{co}), i.e., $\disp\la \psi(t), \ou(t)\ra= \max_{u\in[-2,2]}\la \psi(t), u\ra$ for a.e. $t\in[0,\oT].$

{\bf 3.} 
$q(t)=\begin{pmatrix}
q_1(t)\\
q_2(t)
\end{pmatrix}=\begin{pmatrix}
p_1\\
p_2
\end{pmatrix}-\disp\int_{(t,\oT]}d\gamma(\tau)\begin{pmatrix}
\frac{1}{\sqrt{2}}\\
\frac{1}{\sqrt{2}}
\end{pmatrix}$ (by $(\eqref{c:6.9})$.

{\bf 4.}
$\nn\begin{array}{lll} 
\begin{pmatrix}
-p_1\\
-p_2
\end{pmatrix}-\eta(\oT)\begin{pmatrix}
\frac{1}{\sqrt{2}}\\
\frac{1}{\sqrt{2}}
\end{pmatrix}=\mu\begin{pmatrix}
\ox_1(\oT)-\alpha\\
0
\end{pmatrix}+\begin{pmatrix}
0\\
\lambda
\end{pmatrix}\;{\rm for\;some\;\mu\ge 0\;{\rm and}\;\lambda\in\R}.\\ 
\eta(\oT)>0 \Longrightarrow \la (\ox_1(\oT),1),(1, 1)\ra= 1-\oT. \\
\disp\dfrac{1}{\oT}\int^{\oT}_{0}\langle p,\dot{\ox}(t)\rangle dt= \mu \;\;{\rm(by \eqref{42})}.
\end{array}
\right.$

{\bf 5.}
$\eta(t)=0$\;\mbox{ for a.e. }\;$t\in[0,\oT]$ such that $\ox_1(t)+\ox_2(t)< 1-t$, and 
$\eta(t)>0\Longrightarrow\Big\la q(t),\Big(\begin{matrix}
\frac{1}{\sqrt{2}}\\
\frac{1}{\sqrt{2}}
\end{matrix}\Big)\Big\ra=0$\;\\\textrm{ for a.e. }\;$t\in[0,\oT]$ (by \eqref{41}).

{\bf 6.} $(\mu, p,\n\gg_0\en_{TV},\n\gg_>\en_{TV})\neq 0$ and $\gg=\gg_>+\gg_0$.

{\bf 7.} The following support conditions hold:
\begin{equation*}
\left\{\begin{array}{ll} 
\mathrm{supp}(\gg_0)\cup\;\mathrm{supp}(\gg_>)\subset
\{t\;|\:\ox_1(t)+ \ox_2(t)=1-t\}\\
\mathrm{supp}(\gg_>)\cap \mathrm{int}(E_0)=\emptyset,\;\textrm{ provided that }\;\mathrm{int}(E_0)==\mathrm{int}(\{t\;|\;\eta(t)>0\})\;\textrm{ is nonempty}. 
\end{array}
\right.
\end{equation*}
To utilize these conditions, first we obtain from {\bf 1} that
$\dot{p}(t)=0$ for a.e. $t\in[0,\oT]$, which implies that $p=
\begin{pmatrix}
p_1\\
p_2
\end{pmatrix}$ is constant on $[0,\oT]$ and that $\psi(t)=q_2(t)$. It follows from {\bf 2} that $q_2(t)=\psi(t)\in N(\ou(t);U)$, i.e., 
$$q_2(t)\ou(t)=\max_{u\in U}q_2(t)u\;\mbox{ for a.e. }\;t\in [0,\oT].$$ 
\begin{figure}[h]
\centering
\includegraphics[width=2.9in]{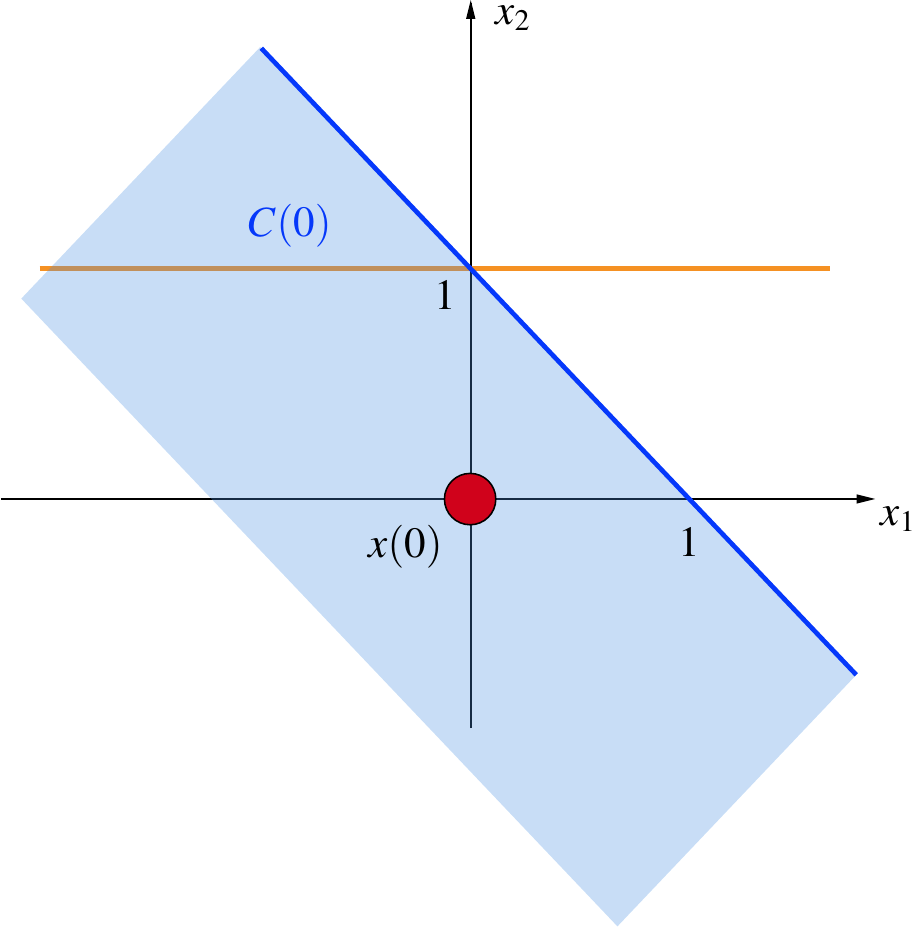}
\caption{The initial point $(0,0)$ belongs to the interior of the moving set $C(t)$}
\label{fig1}
\end{figure}
Note that the initial point $(0,0)$ belongs to the interior of the moving set $C(t)$ and the boundary of $C(0)$ intersects the $y$-axis exactly at the target line $\{y=1\}$; see Figure~\ref{fig1}. Since the unconstrained movement of $x$ may occur only along the $y$-axis and $\mathrm{bd}C(t)\cap \{(0,y)\;|\;y\in \R\}=\{(0,1-t)\}$, it follows that at the final time $\oT$ we have $(\ox_1(\oT),\ox_2(\oT))\in \mathrm{bd}\,C(\oT)$, i.e., $\ox_1(\oT)=-\oT$. Denote by $t_h\in(0,\oT]$ the first time when $\ox(t)$ hits the boundary $\mathrm{bd}\,C(t)$; see Figure~\ref{fig2}.\vspace*{-0.05in}

Observe that $C(t)=C_0-v(t)$, where $C_0=\{(x_1,x_2)\;|\;x_1+x_2\leq 1\}$ and  $v(t)=(\frac{t}{2},\frac{t}{2})$. With the change of variables $y:=x+v$, we get that
$y(t)\in C_0$ for all $t$, $ N(y;C_0)=N(x;C(t))$, and 
$$
\dot y(t)\in -N(y(t);C_0)+\dot v(t)+(0,u(t))=-N(y(t);C_0)+\Big(\dfrac{1}{2}, \dfrac{1}{2}+u(t)\Big).$$
It is well known (see, e.g., \cite[Theorem~10.1.1]{jpa}) that in this case, the dynamics coincides with 
\begin{equation}
\label{A}
\dot y(t) = \proj_{T(y(t);C_0)}\Big(\dfrac{1}{2}, \dfrac{1}{2}+u(t)\Big).
\end{equation}
This implies, in particular, that the normal cone is {\em active} at a.e. $t$ with $u(t)>-1$ and $x_1(t)+x_2(t)=1-t$. In other words, for such $t$ the sweeping dynamics reads as 
$$
\begin{aligned}
\begin{pmatrix}
\dot x_1\\
\dot x_2
\end{pmatrix}=\begin{pmatrix}
\dot y_1\\
\dot y_2
\end{pmatrix}-\dot v= \begin{pmatrix}
\dot y_1\\
\dot y_2
\end{pmatrix}-\dfrac{1}{2}\begin{pmatrix}
1\\
1
\end{pmatrix}&=\begin{pmatrix}
\frac{1}{2}\\[1ex]
\frac{1}{2}+u
\end{pmatrix}-\left\langle\begin{pmatrix}
\frac{1}{2}\\[1ex]
\frac{1}{2}+u
\end{pmatrix},
\begin{pmatrix}
\frac{1}{\sqrt{2}}\\[1ex]
\frac{1}{\sqrt{2}}
\end{pmatrix}\right\rangle
\begin{pmatrix}
\frac{1}{\sqrt{2}}\\[1ex]
\frac{1}{\sqrt{2}}
\end{pmatrix}-\begin{pmatrix}
\frac{1}{2}\\[1ex]
\frac{1}{2}
\end{pmatrix}\\
& =\begin{pmatrix}
0\\
u
\end{pmatrix}-\begin{pmatrix}
\frac{1}{2}\\[1ex]
\frac{1}{2}
\end{pmatrix}-\begin{pmatrix}
\frac{u}{2}\\[1ex]
\frac{u}{2}
\end{pmatrix}=\begin{pmatrix}
-\frac{1+u}{2}\\[1ex]
\frac{u-1}{2}
\end{pmatrix}
\end{aligned}
$$
telling us that $\dot x_1(t)=-\frac{1+u}{2}$ and $\dot x_2(t)= \frac{u-1}{2}$ provided that the normal cone constraint is active on the whole interval $[t_h,t]$. Since $x_1(t_h)=0$, it follows from the above that
$$
x_2(t)=-x_1(t)+1-t= \dfrac{t-t_h}{2}+\dfrac{1}{2}\int^t_{t_h} u(s)ds +1-t= 1-\dfrac{t+t_h}{2}+\dfrac{1}{2}\int^t_{t_h} u(s)ds.
$$
This implies by ${\bf 5}$ that for a.e. $t$ with $x_1(t)+x_2(t)=1-t$ and $u>-1$ a.e. on $[t_h,t]$ we have
\begin{equation*}
\eta(t)=\left\langle-\(\dot x(t)-(0,u(t)\), \begin{pmatrix}
\frac{1}{\sqrt{2}}\\
\frac{1}{\sqrt{2}}
\end{pmatrix}\right\rangle=\left\langle\begin{pmatrix}
\frac{u+1}{2}\\
\frac{u+1}{2}
\end{pmatrix}, 
\begin{pmatrix}
\frac{1}{\sqrt{2}}\\
\frac{1}{\sqrt{2}}
\end{pmatrix}\right\rangle=\dfrac{u+1}{\sqrt{2}}.
\end{equation*}
Observe also that, in the limit case where $\ou(t)\equiv -1$, the trajectory moves downwards parallel to the $y$-axis while always remaining in $\mathrm{bd}\,C(t)$ without any action of the normal cone. Moreover, it follows from \eqref{A} that if $\ou(t)\leq -1$ on an interval $I$, then the normal cone constraint doesn't play any role, and on $I$ the sweeping dynamics reduces to the classical one $\dot x=g(x,u)$.\vspace*{-0.05in}

Observe further that the cost function depends only on the final time and on the final position. Without loss of generality, we assume for simplicity of computations that the control is constant as long as the dynamics remains the same. In particular, we may suppose that the control is constant, say $\ou\equiv u_1$, on $[0,t_h]$. This implies that $u_1>-1$, otherwise $t_h=\infty$. Then the final time $t_h$ can be expressed through the formula
\begin{equation}
\label{th}
t_h=\dfrac{1}{u_1+1}, 
\end{equation}
and if $\ox_1(\oT)+\ox_2(\oT)\in \mathrm{bd}\,C(\oT)$, i.e., $\ox_1(\oT)+\ox_2(\oT)=1-\oT$, we deduce from $\ox_2(\oT)=1$ that 
\begin{equation}
\label{x1T}
\ox_1(\oT)=-\oT.
\end{equation}
\begin{figure}[ht]
    \centering
\includegraphics[width=2.9in]{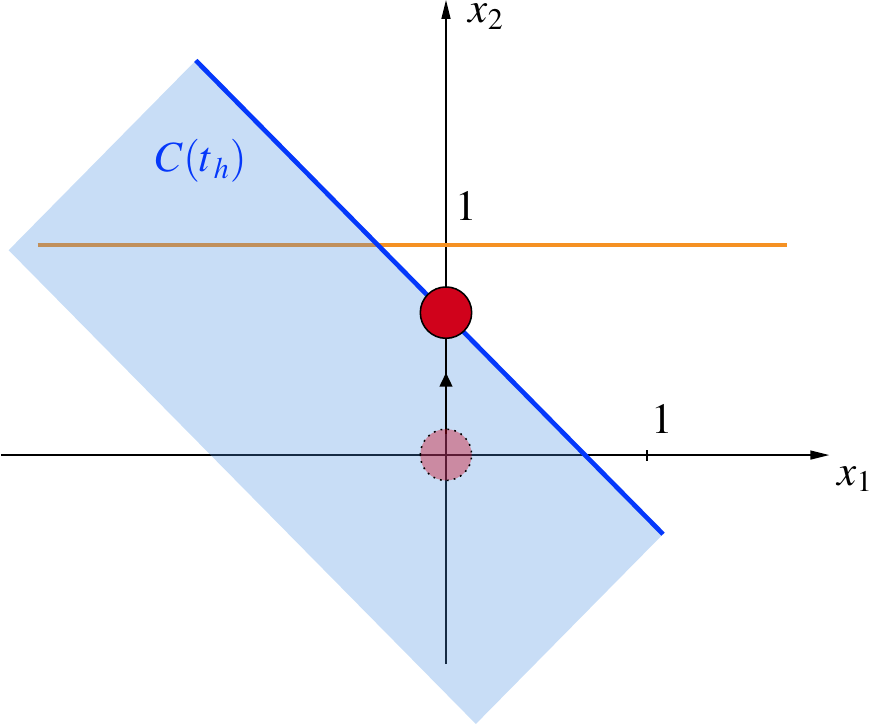}
    \caption{The first time when $x(t)$ hits $\mathrm{bd}\,C(t)$}
    \label{fig2}
\end{figure}

There are the following three cases to consider:

{\bf Case~1:} $t_h=\oT$, i.e., $\ox_1(t)\equiv 0$ and $\ox_2(t)<1$ for all $t<\oT$, but $\ox_2(\oT)=1$. This is impossible since the boundary of $C(t)$ is moving down and $(0,1)\in \mathrm{bd}\,C(0)$.

{\bf Case~2:} $t_h<\oT$ and $\ou(t)>-1$ for a.e. $t\in [t_h,\oT]$. In this case, as previously discussed, the normal cone is active on the whole interval $[t_h,\oT]$. Therefore, by {\bf 5} we have $q_1(t)+q_2(t)=0$ a.e.\ on $[t_h,\oT]$ giving us $p_1+p_2=\sqrt{2}\gg((t,\oT])=\sqrt{2}\gg_0((t_h,\oT])$ for a.e. $t\in [t_h,\oT]$. This implies, in particular, that the measure $\gg$ is concentrated at the final time $\oT$. Thus we get:

{\bf(a)} $p_1+p_2=\sqrt{2}\gg_0(\{\oT\})$.

{\bf(b)} $p_1\ox_1(\oT)+p_2=\mu\oT$ obtained by {\bf 4}. 

{\bf(c)} $p_2=\dfrac{-\eta(\oT)}{\sqrt{2}}-\lambda$ obtained by {\bf 4}.

{\bf(d)} $-p_1+p_2=\mu(\ox_1(\oT)-\al)-\lm$ obtained by {\bf 4} and {\bf (c)}.\\[1ex]
Assuming first that $\mu=0$ yields 
$$p_1=p_2+\lambda,\;\mbox{ and }\;p_2(\ox_1(\oT)+1)=-\lambda\ox_1(\oT).$$
By the nontriviality condition  
and (b), we have  that $p_2\neq 0$, and so $\ox_1(\oT)=-1$ with $\oT=1$ provided that $\lambda=0$; see Figure~\ref{fig3}. In this case, the cost is $\varphi(-1,1)=1+\frac{1}{2}(-1-\alpha)$.
\begin{figure}[h]
\centering
\includegraphics[width=3.1in]{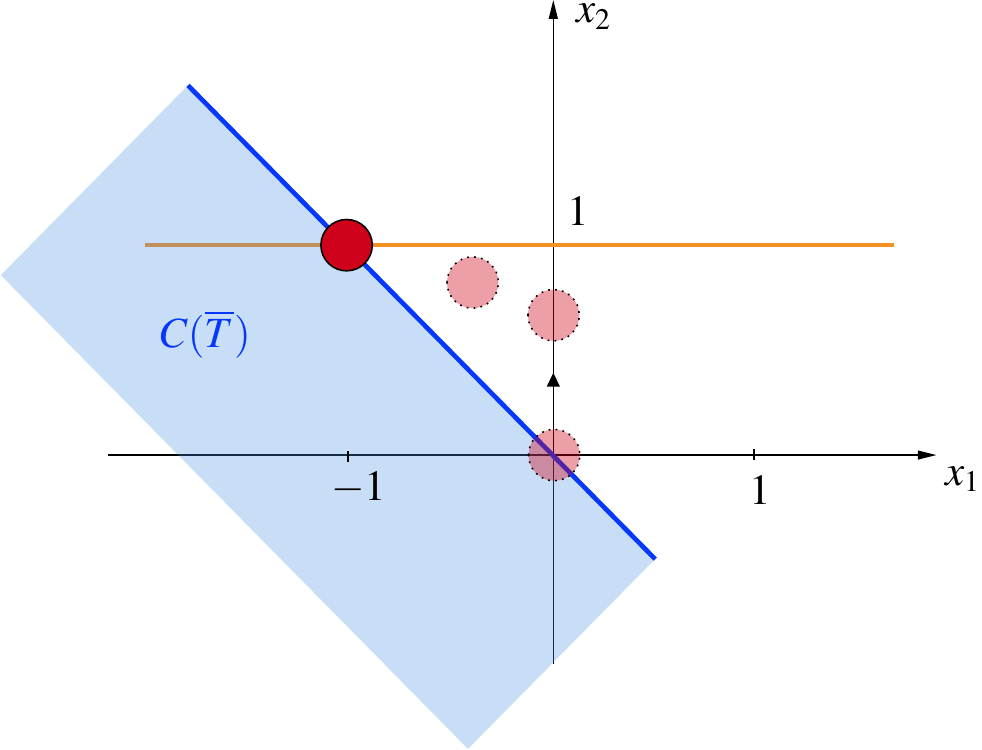}
\caption{When $t_h<\oT$ and $\ou(t)>-1$ for a.e. $t\in [t_h,\oT]$}
\label{fig3}
\end{figure}
Moreover, for any $\lm$, our computation of 
$\eta (t)$ yields $p_2=-\frac{\eta(\oT)}{\sqrt{2}}-\lambda$ 
by (c). 
Involving (a), we deduce that 
$\gg_0(\{\oT\})=\sqrt{2}p_2+\frac{\lambda}
{\sqrt{2}}$. 
Thus $q_2=p_2-\frac{1}{\sqrt{2}} \gg_0(\{\oT\})=$
$-\frac{\lambda}{2}$, 
which yields $u_1=u_2=2$, since $u_2>-1$ provided that $\lambda\neq 0$. 
If instead $\mu=1$, then items (a)--(d) bring us to the following system:
\begin{equation*}
\left\{\begin{array}{lll} 
p_1+p_2=\sqrt{2}\gg_0(\{\oT\}),\\
\oT= p_1\ox_1(\oT)+p_2=-\oT p_1+p_2,\\
p_2=-\dfrac{\eta(\oT)}{\sqrt{2}}-\lambda=-\dfrac{1+u_2}{2}-\lambda,\\
-p_1+p_2=\ox_1(\oT)-\al-\lm.
\end{array}
\right.
\end{equation*}
We get therefore that $\oT(1+p_1)=p_2$. The two cases may happen: either $p_1=-1$, or $p_1\neq -1$ and $\oT=\frac{p_2}{1+p_1}$. Then $q_2=p_2-\frac{1}{\sqrt{2}}\gg_0(\{\oT\})=p_2-(p_1+p_2)=-p_1$. If $p_1\neq 0$ then $u_1=u_2=2$, and in this case there is only one trajectory that satisfies the necessary conditions, with the final cost 
\begin{equation}
\label{C1}
\varphi(\ox_1(\oT),\oT)=1+\dfrac{1}{2}(-1-\al)^2,
\end{equation}
that agrees with the cost that we computed if $\mu=0$. If $p_1=0$, then $p_2=\oT$, so that $\oT=-\frac{\alpha+\lm}{2}$, and $q_2=p_2-\frac{1}{\sqrt{2}}\gg_0(\{\oT\})=p_2-\frac{1}{\sqrt{2}}\frac{p_2}{\sqrt{2}}=\frac{p_2}{2}\neq 0$, which yields $u_1=u_2=2$. It follows that $\lm=-p_2-\frac{1+u_2}{2}=-p_2+\frac{3}{2}=-\oT+\frac{3}{2}$, and so $\oT=-\al-\frac{3}{2}$. Since in this case $(\ox_1(\oT),\ox_2(\oT))\in \mathrm{bd}\,C(\oT)$, we know that $\ox_1(\oT)=-\oT$. The shortest time to reach the target with the strategy of case 2 is obviously $\oT=1$, so that $\oT=-\al-\frac{3}{2}$ makes sense only if $\al\leq -\frac{5}{2}$. The cost is 
\begin{equation}\label{6.44b}
\varphi\Big(\al+\dfrac{3}{2}, -\al-\dfrac{3}{2}\Big)=-\al-\dfrac{3}{8}
\end{equation}
that coincides with the cost in \eqref{C1} for $\al= -\frac{5}{2}$, while it performs better if $\al<  -\frac{5}{2}$.

{\bf Case~3:} The third possibility occurs when the normal cone is active on some interval $[t_h,\tau]$, while it is not active on the nontrivial interval $(\tau,\oT]$. By the previous  computation, the latter means that $u_(t)\leq -1$ on $(\tau,\oT]$ implying that $\dot \ox_1(t)=0$ and $\dot\ox_2(t)=u_3\in [-2,-1]$ on $(\tau, \oT)$. In this case,  we have 
\begin{equation}
\label{x1}
\ox_1(\oT)=\ox_1(\tau)= -\dfrac{1+u_2}{2}(\tau-t_h)=-\dfrac{1+u_2}{2}\Big(\tau-\dfrac{1}{u_1+1}\Big),
\end{equation}
together with $\ox_1(\oT)=-\oT$ if $u_3\equiv-1$, and 
\begin{equation}
\label{x2}
\begin{aligned}
1=\ox_2(\oT)=\ox_2(\oT)-\ox_2(\tau)+\ox_2(\tau)&= (\oT-\tau)u_3-\ox_1(\tau)+1-\tau\\
&=(\oT-\tau)u_3+\dfrac{1+u_2}{2}\Big(\tau-\dfrac{1}{u_1+1}\Big)+1-\tau.
\end{aligned}
\end{equation}
Since $\eta(\oT)=0$, the necessary conditions {\bf 4} and {\bf 5} can be rewritten, respectively, as

{\bf 4'}. $\nn\begin{array}{lll} 
-(p_1,p_2)=\big(\mu(\ox_1(\oT)-\al\big), \lambda),\\
p_1\ox_1(\oT)+p_2=\mu\oT.
\end{array}
\right.$ 

{\bf 5'}. $q_1(t)+q_2(t)=0$ a.e.\ on $[t_h,\tau)$.\\[1ex]
Moreover, we have $\gg=\gg_0$ on $[t_h,\tau)$, $\gg(\{\tau\})=\gg_0(\{\tau\})+\gg_>(\{\tau\}),$ and $\gg=\gg_>$ on $(\tau,\oT]$, where $\gg_>=0$ if $u_3>-1$. It follows therefore that
$$p_2=-\lm,\;p_1=\mu(\al-\ox_1(\tau)),\;\mbox{ and }\;p_1\ox_1(\tau)=\mu\oT+\lm.$$
Furthermore, by {\bf 5'} we get $0=p_1+p_2-\sqrt{2}\big(\gg_0((t,\tau])+\gg_>([\tau,\oT])\big)$ a.e.\ on $[t_h,\tau)$, which implies that the support of $\gg_0$ is contained in $\{\tau\}$, telling us that
$$
p_1+p_2=\sqrt{2}\big(\gg_0(\{\tau\})+\gg_>([\tau,\oT])\big)\;\mbox{ and }\;q_2=p_2-\frac{1}{\sqrt{2}}\Big(\gg_0(\{\tau\})+\gg_>([\tau,\oT])\big)\;\mbox{on }[0,\tau].$$ 
Since $\ox_1(\oT)< 0$ from the structure of the problem, it follows from {\bf 4'} that $p_1\neq 0$, and so $\gg_0(\{\tau\})+\gg_>([\tau,\oT]))<0$ provided that $\mu=1$ and $\lm=0$.
The maximization condition {\bf 2} now reads as
$$
0>q_2\in N(\ou(t);U)\;\mbox{ for a.e. }\;t\in (0,\tau),
$$
which yields $\ou(t)\equiv 2$ on $(0,\tau)$ when $\mu=1$ and $\lm=0$. Thus $t_h=\frac{1}{3}$ by \eqref{th}. Taking $p_2=-\lm\neq 0$ gives us necessarily  that $\mu\neq 0$. Then if $p_1=0$, then $\ox_1(\oT)=\al$ and the cost is $\varphi(\al,-\al)=-\al$, while for $p_1\neq 0$ we get $q_2\neq 0$ and $\ou=2$ on $(0,2)$ as previously discussed. Continuing our analysis, consider the cases:

{\bf(i)} If $\mu =0$, then $p_1=p_2=0$, and for $t\in (t_h,\tau)$ we have $0=q_1(t)+q_2(t)=-\sqrt{2}\big(\gg_0(\{\tau\})+\gg_>([\tau,\oT])\big)$ implying that $q_1=q_2=0$. There are two possibilities: either $u_3=-1$, or $\gg_>=0$, that implies $\gg_0=0$, which violates the nontriviality condition.

{\bf(ii)} If $\mu =1$, then the two subcases may occur on the final interval $(\tau,\oT)$:\\[1ex]
{\bf(a)} If $u_3=-1$, then $\eta\equiv 0$ on $[\tau,\oT]$, while the trajectory keeps contact with the boundary of the moving constraint. In this case $\ox_1(\oT)=-\oT$ as already explained,. Taking into account the expressions of $x_1$ and $x_2$, we obtain that $\oT=\frac{3}{2}\tau-\frac{1}{2}$, and the final cost is 
\begin{equation*}
\varphi(\ox_1(\oT),\oT)=\frac{3}{2}\tau-\frac{1}{2}+\frac{1}{2}\(\frac{3}{2}\tau-\frac{1}{2}+\al\)^2,      
\end{equation*}
which attains its minimum at $\tau=-\frac{1}{3}-\frac{2}{3}\al$. In this case, the optimal cost is 
\begin{equation}
\label{C2}
J=-\frac{1}{2}-\al.
\end{equation}
This solution is meaningful if and only if $\al<-2$ since in this case $\oT>1$ as previously discussed. Comparing \eqref{C1} and \eqref{C2} shows that \eqref{C2} performs better than \eqref{C1}.\\[1ex]
{\bf(b)} If $u_3<-1$, then the dynamics is $\dot \ox_1=0$, $\dot \ox_2=u_3$. In this case, $p_2=q_2=0$ on $(\tau, \oT)$, and so no information comes from the maximization condition. However, since the second summand of the final cost depends only on $\ox_1(\tau)$, i.e., on the second switching time $\tau$, and we already proved that in the two first intervals the optimal (constant) control corresponds to the maximal speed, it is obvious that the best performance is obtained when $u_3\equiv -2$. Hence the strategy with $u_3=-1$, whose cost was computed in \eqref{C2}, is not optimal. Therefore, in the case where $u_3\equiv -2$, the optimal performance can be computed directly. We deduce from \eqref{x1} and \eqref{x2}, with $u_1=u_2=2$ and $u_3=-2$, that
\begin{equation*}
\nn\begin{array}{lll} 
\ox_1(\oT)=\ox_1(\tau)=-\dfrac{3}{2}\Big(\tau-\dfrac{1}{3}\Big),\\
\oT=\dfrac{5}{4}\tau-\dfrac{1}{4},\\
\ox_2(\tau)=\dfrac{\tau+1}{2}.\\
\end{array}
\right.
\end{equation*}
\begin{figure}[htbp]
\centering
\begin{subfigure}[b]{0.47\textwidth}
\includegraphics[width=\textwidth]{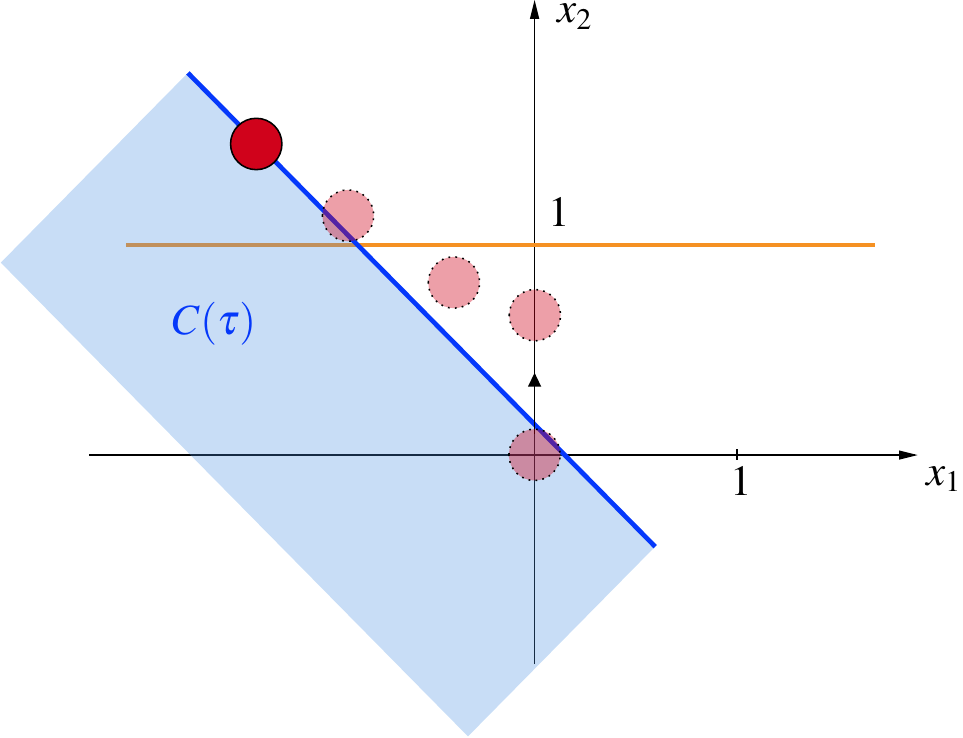}
\caption{The object is at the second switching time $\tau$}
\end{subfigure}
\hfill
\begin{subfigure}[b]{0.49\textwidth}
\includegraphics[width=\textwidth]{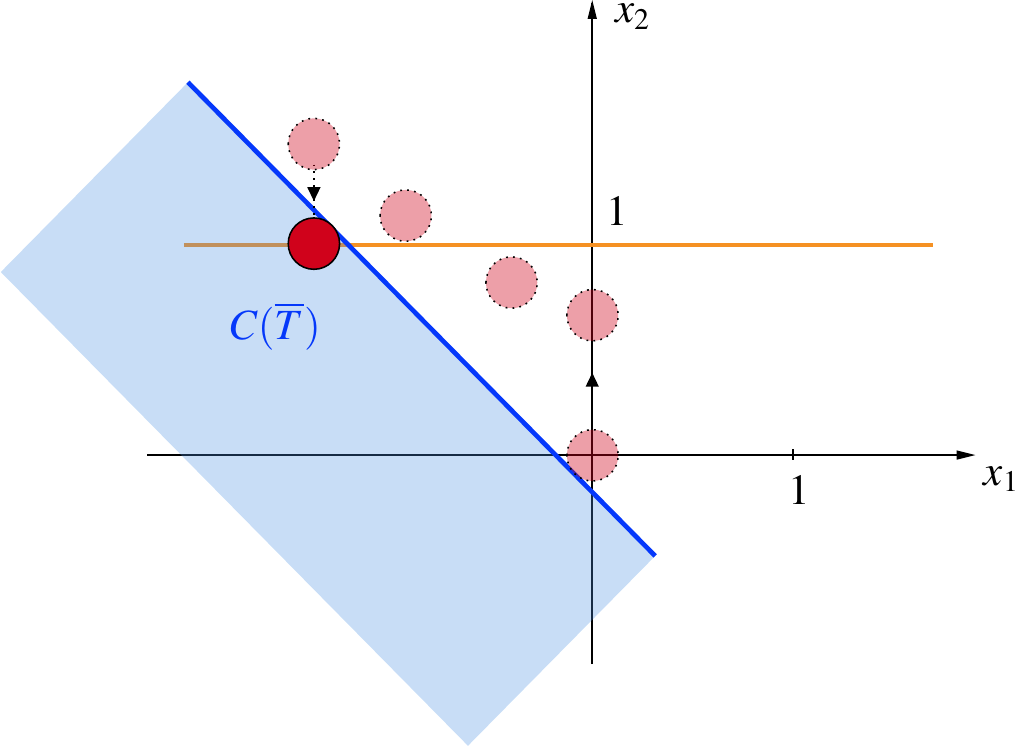}
\caption{The object reaches the target}
\end{subfigure}

\caption{ The strategy with a downwards switching}
\label{fig4}
\end{figure}
Observe that this strategy makes sense if and only if $\ox_2(\tau)>1$, which is equivalent to $\oT>\tau$. The final cost can be computed accordingly as a function of $\tau$ and of the parameter $\al$:
$$
\varphi(\ox_1(\oT),\oT)=\dfrac{5}{4}\tau-\dfrac{1}{4}+\dfrac{1}{2}\Big(\dfrac{3}{2}\tau-\dfrac{1}{2}+\al\Big)^2.
$$
Thus the optimal switching time $\bar{\tau}$ is $$\bar{\tau}=-\dfrac{2}{9}-\dfrac{2}{3}\al,$$
which makes sense if and only if $\bar{\tau}>1$, i.e. $\al<\frac{-11}{6}$, and the optimal value is expressed by the formula
\begin{equation}
\label{C3}
\varphi(\ox_1(\bar \tau),\bar \tau )=-\dfrac{13}{72}-\dfrac{5}{6}\al.  
\end{equation}
The cost value in \eqref{C3} is strictly smaller than the cost values in \eqref{C1}, \eqref{C2}, and \eqref{6.44b} by $\al<-2$, i.e., the strategy with a downwards switching perform better than the strategy without switching for the same parameter value $\al$; see Figure~\ref{fig4}.
\end{example}\vspace*{-0.15in}

\section{Concluding Remarks}\label{sec:Conclusions}\vspace*{-0.05in}

This paper develops a novel version of the method of discrete approximations to study a challenging class of free-time optimal control problems for sweeping processes with the nonsmooth controlled dynamics and endpoint constraints. Being of its own interest, this method combined with powerful tools of (first- and second-order) variational analysis and generalized differentiation allows us to derive a new set of necessary optimality conditions for discrete-time and continuous-times sweeping control systems governed by moving polyhedra. The obtained illustrated by an instructive example. Their applications to practical models of marine surface vehicles and nanoparticle dynamics will be presented in the forthcoming paper.
\end{document}